\documentclass[12pt]{article}

\usepackage{subcaption}

\usepackage{float}

\usepackage{multirow}

\usepackage{tikz}

\usetikzlibrary{trees}

\usetikzlibrary{calc,backgrounds,arrows,matrix}
\usetikzlibrary{shapes.geometric} 
\usetikzlibrary{shapes,snakes}

\usepackage{enumerate}

\usepackage{color}
\definecolor{lightblue}{rgb}{0,0.2,0.5}

\usepackage[colorlinks=true, urlcolor=blue,linkcolor=blue, citecolor=lightblue]{hyperref}

\usepackage[round,sort,comma,numbers]{natbib}

\bibpunct{\textcolor{lightblue}{(}}{\textcolor{lightblue}{)}}{,}{a}{}{;}

\usepackage{amssymb,amsmath}

\usepackage{graphicx}

\usepackage{graphics}
\usepackage{color}

\DeclareMathAlphabet{\eufrak}{U}{}{}{}
\SetMathAlphabet\eufrak{normal}{U}{euf}{m}{n}
\SetMathAlphabet\eufrak{bold}{U}{euf}{b}{n}

\allowdisplaybreaks

 \def\qu{{\mathord{\mathbb Z}}}

 \def\T{{\mathrm{{\rm T}}}}

 \def\inte{{\mathord{\mathbb R}}}

 \def\inte{{\mathord{\mathbb N}}}

 \def\sZZ{{\rm Z\kern-.45em{}Z}}

 \def\sQQ{{\kern 0.27em \vrule height1.45ex width0.03em depth0em
           \kern-0.30em \rm Q}}
 \def\qu{{\mathchoice
         {\sQQ}
         {\sQQ}
   {\kern 0.225em \vrule height1.05ex width0.025em depth0em \kern-0.25em \rm Q}
   {\kern 0.180em \vrule height0.78ex width0.020em depth0em \kern-0.20em \rm Q}
         }}
 \def\sGG{{\kern 0.27em \vrule height1.45ex width0.03em depth0em
           \kern-0.30em \rm G}}
 \def\gg{{\mathchoice
         {\sGG}
         {\sGG}
   {\kern 0.225em \vrule height1.05ex width0.025em depth0em \kern-0.25em \rm G}
   {\kern 0.180em \vrule height0.78ex width0.020em depth0em \kern-0.20em \rm G}
         }}

 \newtheorem{prop}{Proposition}[section]
 \newtheorem{lemma}[prop]{Lemma}
 \newtheorem{definition}[prop]{Definition}
 
 \newtheorem{theorem}[prop]{Theorem}

\numberwithin{equation}{section}

 \newcounter{hyp}
\newenvironment{Proofy}{\removelastskip\par\medskip \noindent{\em Proof of Theorem} \rm}{\penalty-20\null\hfill$\square$\par\medbreak}

\newenvironment{Proofz}{\removelastskip\par\medskip \noindent{\em Proof \hskip-0.3cm } \rm}{\penalty-20\null\hfill$\square$\par\medbreak}

\newenvironment{Proof}{\removelastskip\par\medskip \noindent{\em Proof.} \rm}{\penalty-20\null\hfill$\square$\par\medbreak}

\def\bprf{\begin{Proof}}
\def\nprf{\end{Proof}}
\def\bdes{\begin{description}}
\def\ndes{\end{description}}

\newcommand{\E}{\mathbb{E}}

\newcommand{\real}{\mathbb{R}}

\def\og{\leavevmode\raise.3ex
     \hbox{$\scriptscriptstyle\langle\!\langle$~}}
\def\fg{\leavevmode\raise.3ex
     \hbox{~$\!\scriptscriptstyle\,\rangle\!\rangle$}~}

\title{\Huge
Semilinear fractional elliptic PDEs with gradient nonlinearities on open balls: existence of solutions and probabilistic representation 
}

\author{
  Guillaume Penent\footnote{\href{mailto:PENE0001@e.ntu.edu.sg}{pene0001@e.ntu.edu.sg}}
  \qquad 
      Nicolas Privault\footnote{
\href{mailto:nprivault@ntu.edu.sg}{nprivault@ntu.edu.sg}
}
  \\
\small
Division of Mathematical Sciences
\\
\small
School of Physical and Mathematical Sciences
\\
\small
Nanyang Technological University
\\
\small
21 Nanyang Link, Singapore 637371
}

\allowdisplaybreaks

\usepackage{empheq}

\usepackage{bbm}

\makeatletter
\newcommand*\rel@kern[1]{\kern#1\dimexpr\macc@kerna}
\newcommand*\widebar[1]{
  \begingroup
  \def\mathaccent##1##2{
    \rel@kern{0.8}
    \overline{\rel@kern{-0.8}\macc@nucleus\rel@kern{0.2}}
    \rel@kern{-0.2}
  }
  \macc@depth\@ne
  \let\math@bgroup\@empty \let\math@egroup\macc@set@skewchar
  \mathsurround\z@ \frozen@everymath{\mathgroup\macc@group\relax}
  \macc@set@skewchar\relax
  \let\mathaccentV\macc@nested@a
  \macc@nested@a\relax111{#1}
  \endgroup
}
\makeatother

\let\oldcitet=\citet
\let\oldcitep=\citep
\renewcommand{\cite}[1]{\textcolor[rgb]{0,0,1}{\oldcitet{#1}}}
\renewcommand{\citet}[1]{\textcolor[rgb]{0,0,1}{\oldcitet{#1}}}
\renewcommand{\citep}[1]{\textcolor[rgb]{0,0,1}{\oldcitep{#1}}}

\renewcommand{\leq}{\leqslant}

\renewcommand{\geq}{\geqslant}

\begin{document}

\maketitle

\baselineskip0.6cm

\vspace{-0.6cm}

\begin{abstract}
 We provide sufficient conditions for the existence of viscosity solutions of fractional semilinear elliptic PDEs of index $\alpha \in (1,2)$ with polynomial gradient nonlinearities on $d$-dimensional balls, $d\geq 2$. Our approach uses a tree-based probabilistic representation of solutions and their partial derivatives using $\alpha$-stable branching processes, and allows us to take into account gradient nonlinearities not covered by deterministic finite difference methods so far. In comparison with the existing literature on the regularity of solutions, no polynomial order condition is imposed on gradient nonlinearities. Numerical illustrations demonstrate the accuracy of the method in dimension $d=10$, solving a challenge encountered with the use of deterministic finite difference methods in high-dimensional settings. 
\end{abstract}

\noindent
    {\em Keywords}:
Elliptic PDEs, 
semilinear PDEs,
fractional Laplacian,
gradient nonlinearities, 
stable processes,
branching processes,
Monte-Carlo method.

\noindent
    {\em Mathematics Subject Classification (2020):}
35J15, 
35J25, 
35J60, 
35J61, 
35R11, 
 35B65, 
 60J85, 
 60G51, 
 60G52, 
 65C05, 
 33C05. 
        
\baselineskip0.7cm

\parskip-0.1cm

\section{Introduction}
The study of solutions of nonlocal and fractional elliptic
partial differential equations (PDEs) is an active research topic
which has attracted significant attention over the past decades.
In the case of the classical (local) Laplacian, viscosity solutions
 of fully nonlinear second-order elliptic PDEs
 have been constructed in \cite{ishii1989}
 by the Perron method.

\medskip

 On the other hand, nonlocal elliptic PDEs
 can be solved using weak solutions, see Definition~2.1 in \cite{ros-oton},
 or viscosity solutions, see \cite{servadei2014}
 and Remark~2.11 in \cite{ros-oton}. 
 Weak solutions can be obtained from the Riesz representation 
 or Lax-Milgram theorems as in \cite{felsinger}, \cite{ros-oton2016}. 
 See also \cite{barles2}
 for the use of the Perron method,
 and \cite{bony} for semi-group methods
 applied to second-order elliptic integro-differential
 PDEs. 
 
      \medskip

 Given $d\geq 1$, let 
      $$
 \Delta_\alpha  u = - ( - \Delta )^{\alpha /2} u
 = \frac{4^{\alpha /2} \Gamma ( d/2 + \alpha /2 )}{\pi^{d/2}|\Gamma (- \alpha /2 )|}
 \lim_{r\rightarrow 0^+}
 \int_{\mathbb{R}^d \setminus B(x,r)} \frac{u(  \ \! \cdot +z)-u(z)}{|z|^{d+\alpha }}dz,
 $$
  denote the fractional Laplacian
 on $\real^d$ with parameter $\alpha \in (0,2)$,
 see, e.g., \cite{tendef},
      where $\Gamma (p) : = \displaystyle
   \int_0^\infty e^{-\lambda x} \lambda^{p-1} d\lambda$ 
   is the gamma function
   and $|z|$ is the Euclidean norm of $z\in \real^d$. 

   \medskip

   For problems of the form
      $$
      \Delta_\alpha u(x) + f(x) =0,
      $$
 with $u=\phi$ on $\real^d \setminus D$, 
 where $D$ is an open bounded domain in $\real^d$,
 existence of viscosity solutions 
   has been proved in \cite{servadei2014}
   under smoothness assumptions on 
   $f,\phi$,
   see also \cite{felsinger}, resp. \cite{mou}, for the existence of
   weak solutions, resp. viscosity solutions, with nonlocal operators. 
   Regarding problems of the form
      $$
      \Delta_\alpha u(x) + f(x,u(x)) =0,
      $$
      existence of non trivial solutions
      with $u=0$ outside an
      open bounded domain $D$ with
      Lipschitz boundary in $\real^d$ 
 has been considered in
 \cite{servadei} using the mountain pass theorem
 when $f$ is a Carath\'eodory function on $D\times \real^d$
 satisfying a polynomial growth condition of order
 $m\in (1,(d+\alpha)/(d-\alpha ) )$. 

\medskip
   
 The regularity of viscosity solutions of
 semilinear elliptic PDEs of the form
\begin{equation}
\label{eq:0}
  \Delta_\alpha u(x)
  - b(x) \Vert \nabla u(x) \Vert_{\real^d}^{\kappa + \tau } 
  - \Vert \nabla u(x) \Vert_{\real^d}^r 
  = 0, \qquad x \in D, 
\end{equation}
 where $D$ is an open domain of $\real^d$,
 $b$ is in the space $C^{0,\tau }$ of $\tau$-H\"older continuous
 functions, $\tau \in (0,1)$,
 has been considered in \S~4.3 in \cite{barles3},
 namely, as a consequence of Theorem~3.1 therein,
 if $b$ is in the space $C^{0,\tau }$ of $\tau$-H\"older continuous
 and $\kappa , r \in (0,2)$, then
 any bounded viscosity solution $u$ of \eqref{eq:0}
 is $\beta$-H\"older continuous
 for small enough $\beta$,
 see also \S~4.1.2 of \cite{barles4} for Lipschitz
 regularity in the case of mixed 
 local and fractional Laplacians. 

 \medskip 

 More recently, the Lipschitz regularity of viscosity solutions of
\begin{equation}
\label{eq:0-1}
  \Delta_\alpha u(x)
  + f(x,\nabla u(x) ) 
  = 0
\end{equation} 
on $D= B(0,R)$
 the open ball
 of radius $R>0$ in $\mathbb{R}^d$,
has been obtained in Theorem~2.1 of \cite{biswap}, 
provided that $f \in {\cal C}(\real^d \times \real^d )$
satisfies a power-type growth
condition of order $m \in (0,\alpha + 1)$ in 
$\nabla u(x)$, while this bound 
can be lifted under an extra coercivity condition on $H$.

\medskip 

 In this paper, we consider the class of semilinear elliptic problems
 on $B(0,R)$ of the form
\begin{equation}
\label{eq:1}
\begin{cases}
  \displaystyle
  \Delta_\alpha  u(x) +
  f(x,u(x) ,\nabla u(x) )
  = u(x), \qquad x \in B(0,R), 
  \medskip
  \\
u(x) = \phi(x), \qquad x \in \mathbb{R}^d \backslash B(0,R),
\end{cases}
\end{equation}
where
\begin{itemize}
\item
  $f(x,y,z)$ is a polynomial nonlinearity on
  $\real^d \times \real \times \real^d$, of the form 
  \begin{equation}
\label{jfklds} 
f(x,y ,z) = \sum\limits_{l= (l_0,\ldots , l_m) \in {\cal L}_m} c_l(x) y^{l_0}
\prod_{i=1}^m (b_i(x) \cdot z)^{l_i},
\end{equation}
where
${\cal L}_m$ is a finite subset of $\mathbb{N}^{m+1}$
for some $m\geq 0$, and 
$(c_l(x))_{l= (l_0,\ldots , l_m)\in {\cal L}_m}$, $(b_i(x))_{i=1,\ldots , m}$ 
are bounded continuous functions of $x \in \real^d$,
with $x\cdot z := x_1z_1+\cdots + x_dz_d$, 
\item
  $\phi : \mathbb{R}^d \rightarrow \mathbb{R}$ is a
  bounded Lipschitz 
 function on $\real^d \setminus B(0,R)$. 
\end{itemize} 
 Using a probabilistic approach, we 
 prove the existence of regular viscosity solutions to \eqref{eq:1}
 under the following conditions. 
 We note that, in comparison to the literature quoted above
 on the regularity of solutions,
 no coercivity or maximum growth order condition in $z$ 
 is imposed on $f(x,y,z)$. 
  
 \medskip

\noindent
{\bf Assumption (\hypertarget{A}{$A$})} 
\begin{enumerate}[1)]
\item The boundary condition $\phi$
   belongs to
   the fractional Sobolev space 
  $$H^\alpha ( \real^d ):=
  \left\{ u \in L^2( \real^d ) \ : \
  \frac{|u(x)-u(y)|}{|x-y|^{d/2+\alpha /2}} \in L^2 (\real^d\times \real^d )\right\}
$$
 and is bounded on $\real^d \setminus B(0,R)$.    
\item {\em
  The coefficients $c_l (x)$,
  $l \in {\cal L}_m$, are uniformly bounded functions, i.e., we have 
}
 \begin{equation}
 \nonumber 
   \Vert c_l\Vert_\infty  := \sup_{x \in B(0,R)} |c_l(x)| < \infty, \qquad l = (l_0,\ldots ,l_m) \in {\cal L}_m.
 \end{equation}
\item 
{\em
The coefficients $b_i(x)$, $i = 0,\ldots,m$, are such that
}
$$
\sup_{x \in B(0,R)}\frac{|b_i(x)|}{R-|x|} < \infty,
\qquad 
 i = 1,\ldots , m. 
$$
\end{enumerate}
Theorem~\ref{t1.0} is the main result of this paper.
It is implied by
 Theorem~\ref{t3}, in which we prove the existence of 
 (continuous) viscosity solutions
 for fractional elliptic problems of the form \eqref{eq:1}. 
\begin{theorem}
 \label{t1.0}
 Let $\alpha \in (1,2)$ and $d\geq 2$.
 Under Assumption~(\hyperlink{A}{$A$}), 
 the semilinear elliptic PDE \eqref{eq:1} admits a viscosity solution
 in ${\cal C}^1(B(0,R)) \cap {\cal C}^0 ( \widebar{B}(0,R))$, 
 provided that $R$ and $\max_{l \in {\cal L}_m} \Vert c_l\Vert_\infty $ are 
 sufficiently small. 
\end{theorem}
 Our method of proof relies on the probabilistic representation of
 PDE solutions using stochastic branching processes, 
 as introduced in \cite{skorohodbranching} and 
 \cite{inw}. 
 Probabilistic representations have been applied to the blow-up and existence of
 solutions for parabolic PDEs in \cite{N-S}, \cite{lm}. 
 They have also been recently extended in \cite{claisse}
 to the treatment polynomial nonlinearities in
 gradient terms in elliptic PDEs with (local) diffusion generators,
 following the approach of \cite{labordere} in the parabolic case.
 In this construction, gradient terms
 are associated to tree branches to which
 a Malliavin integration by parts is applied.
 In \cite{penent}, this approach has been extended to
 the treatment of nonlocal pseudo-differential operators of the form $-\eta(-\Delta /2)$ using random branching trees 
 constructed from a L\'evy subordinator, with application to
 parabolic PDEs with fractional Laplacians
 
\medskip

 The existence of viscosity solutions in Theorem~\ref{t1.0} 
 is obtained through a probabilistic representation of the form
\begin{equation}
  \label{otf} 
u(x) := \mathbb{E} [ \mathcal{H}_\phi (\mathcal{T}_{x,0}) ],
\qquad x\in B(0,R),
\end{equation} 
where $\mathcal{H}_\phi (\mathcal{T}_{x,0})$, see \eqref{djsda}, 
is a functional of a random branching tree
$\mathcal{T}_{x,0}$ started at $x\in \real^d$, 
 and constructed in Section~\ref{s2}. 
 The proof of Theorem~\ref{t1.0} also makes use of existence results
 for nonlinear elliptic PDEs with fractional Laplacians 
 derived in \cite{penent2}, see Theorem~1.2 and
 Proposition~3.5 therein. 

 \medskip

 To prove Theorem~\ref{t1.0}, in Proposition~\ref{t4.1} 
 we construct  for each $i=0,\ldots ,m$ a sufficiently integrable functional 
 $\mathcal{H}_\phi (\mathcal{T}_{x,i})$ 
 of a random tree $\mathcal{T}_{x,i}$ such that the
 probabilistic representations 
 $$
 u (x) =  \mathbb{E}\big[ \mathcal{H}_\phi (\mathcal{T}_{x,0})\big], \quad
 b_i(x) \cdot \nabla u(x) = \mathbb{E}\big[ \mathcal{H}_\phi (\mathcal{T}_{x,i})\big],
 \quad x \in  \real^d, 
 $$ 
 hold 
 for
 the viscosity solutions $u(x)$ of \eqref{eq:1}
 and their gradients $b_i(x)\cdot \nabla u(x)$,
 $x \in B(0,R)$, $i=0,\ldots ,m$, 
 under integrability assumptions on 
 $(\mathcal{H}_\phi (\mathcal{T}_{x,i}))_{x\in B(0,R)}$.
 Then, in Proposition~\ref{t4.3} we show that  
 for any $d\geq 2$ and $p\geq 1$,
 $(\mathcal{H}_\phi (\mathcal{T}_{x,i}))_{x\in B(0,R)}$
 is bounded in $L^p (\Omega )$ 
 uniformly in $x\in B(0,R)$, 
 and therefore uniformly integrable, $i=0,\ldots, m$.
 
 \medskip

 For this, we extend arguments of \cite{claisse}
 from the standard Laplacian $\Delta$ and Brownian motion
 to the fractional Laplacian $\Delta_\alpha : =-(-\Delta)^{\alpha /2}$
 and its associated stable process. 
 There are, however, significant differences from the Brownian case.
 In particular, in the stable setting
 we rely on sharp gradient estimates for fractional Green and Poisson kernel 
 proved in \cite{bogdangradient}, 
 and on integrability results for stable process hitting times,
 see \cite{bogdanbarrier}.
 In addition, the behavior of the negative moments of
 stable processes, see \eqref{inverse},
 requires a more involved treatment of integrability
 in small time when showing the boundedness of 
 $(\mathcal{H}_\phi (\mathcal{T}_{x,i}))_{x\in B(0,R)}$
 in $L^p (\Omega )$, for $p\geq 1$. 
 
\medskip

 In addition, we present a Monte Carlo numerical implementation of
 the probabilistic representation \eqref{otf} 
 on specific examples. 
 In comparison with deterministic finite difference
 methods, see e.g. \S~6.3 of \cite{huang-oberman} for the one-dimensional
 Dirichlet problem,
 our approach allows us to take into account gradient nonlinearities. 
 We also note that our tree-based Monte Carlo implementation 
 applies to high-dimensional problems, see Figure~\ref{fig3}
 in dimension $d=10$,  
 whereas the application of deterministic finite difference
 methods to the fractional Laplacian 
 in higher dimensions is challenging,
 see e.g. \cite{huang-oberman}, page~3082.  

 \medskip

 \indent This paper is organized as follows.
 Section~\ref{s2} presents 
 the description of the branching mechanism,
 following the preliminaries on stable processes and
 kernel introduced in Section~\ref{s1}. 
 In Section~\ref{s4} we state and prove our main existence result,
 Theorem~\ref{t3}, 
 for the probabilistic representation of the solution of \eqref{eq:1}. 
 Section~\ref{s5} presents a Monte Carlo numerical implementation
 of our method on specific examples. 

\section{Preliminaries and notation} 
\label{s1} 
 Before proceeding further, we recall some
 preliminary results on
 fractional Laplacians on the ball $B(0,R)$ in $\real^d$.
\subsubsection*{Poisson and Green kernels} 
Given $(X_t)_{t\geq 0}$ an $\real^d$-valued $\alpha$-stable process,
$\alpha \in (0,2)$, we consider the process
$$
X_{t,x} := x+X_t, \qquad t\in \real_+,
$$
 started at $x\in \real^d$, see e.g. \S~1.3.1 in \cite{applebk2}, 
 and the first hitting time 
\begin{equation}
\nonumber 
  \tau_R (x) := \inf \big\{ t \geq 0,~ X_{t,x} \not\in B(0,R) \big\}
\end{equation} 
of $\real^d\setminus B(0,R)$ by $(X_{t,x})_{t\geq 0}$.
 Note that by the bound (1.4) in \cite{bogdanbarrier} we have 
 $\E [ \tau_R (x) ] < \infty$, and therefore
 $\tau_R (x)$ is almost surely finite for all $x\in B(0,R)$.
 The Green kernel $G_R(x,y)$ satisfies 
 \begin{equation}
   \label{G1} 
 \mathbb{E}\left[
   \int_0^{\tau_R (x)} f(X_{t,x}) dt
   \right]
 = \int_{B(0,R)} G_R(x,y) f(y) dy,
 \quad x\in B(0,R), 
\end{equation}
 for $f$ a nonnegative measurable function on $\real^d$. 
 If $\alpha \in (0,2) \setminus \{d\}$, we also have 
\begin{equation} 
 \nonumber 
  G_R(x,y) = \frac{\kappa_\alpha^d}{|x-y|^{d-\alpha}}
  \int_0^{r_0(x,y)} \frac{t^{\alpha/2-1}}{(1+t)^{d/2}} dt,
  \qquad x, y \in B(0,R), 
\end{equation} 
 see Theorem~3.1 in \cite{bucur}, where
$$
 r_0(x,y) := \frac{(R^2-|x|^2)(R^2-|y|^2)}{R^2|x-y|^2}
 \quad \mbox{and}
  \quad \kappa_\alpha^d := \frac{2^{-\alpha} \Gamma (d/2)}{\pi^{d/2} (\Gamma (\alpha/2))^2}. 
$$
 The Poisson kernel 
 $P_R(x,y)$ of the harmonic measure $\mathbb{P}^x\big( X_{\tau_R (x),x} \in dy\big)$
 satisfies 
\begin{equation}
   \label{P1} 
  \mathbb{E} \big[
 f\big( X_{\tau_R (x)}^x \big) \big] 
   = \int_{\real^d \setminus B(0,R)} P_R(x,y) f(y) dy,
 \quad x\in B(0,R), 
\end{equation} 
 for $f$ a nonnegative measurable function on $\real^d$,
 and is given by 
\begin{equation} 
\nonumber 
  P_R(x,y) = \mathcal{A}(d,-\alpha) \int_{B(0,R)} \frac{G_R(x,z)}{|y-z|^{d+\alpha}} dz
\end{equation}
where
$$
\mathcal{A}(d,-\alpha) := 
\frac{2^\alpha \Gamma ( (d+\alpha)/2 ) }{\pi^{d/2} |\Gamma (-\alpha/2)|}.
$$ 
In particular, when $|x|<R$ and $|y|>R$ we have 
$$
P_R(x,y) =
\frac{{\cal C}(\alpha ,d)}{|x-y|^d} 
      \left( \frac{R^2-|x|^2}{|y|^2-R^2} \right)^{\alpha/2}, 
$$
      with ${\cal C}(\alpha ,d) :=\Gamma (d/2) \pi^{-d/2-1} \sin(\pi \alpha/2)$.
      In addition, we have the bounds 
\begin{equation}
\label{poissonbound}
|\nabla_x P_R(x,y)| \leq (d+\alpha)\frac{P_R(x,y)}{R-|x|},
\quad
x\in B(0,R), \ y \in \real^d \setminus \widebar{B}(0,R), 
\end{equation}
 where $\widebar{B}(0,R)$ denotes 
 closed  ball of radius $R>0$ in $\mathbb{R}^d$, 
 see Lemma~3.1 in \cite{bogdangradient}, and 
\begin{equation}
\label{greenbound}
|\nabla_x G_R(x,y)| \leq d \frac{G_R(x,y)}{\min (
  |x-y|,
  R-|x|)},
\quad x,y \in B(0,R), \ x\neq y, 
\end{equation}
 see Corollary~3.3 in \cite{bogdangradient}. 
\subsubsection*{Moments of stable processes} 
In the sequel 
we will need to estimate the negative moments $\E[|X_t|^{-p}]$ 
of an $\alpha$-stable process  $(X_t)_{t\geq 0}$ 
 represented as the subordinated Brownian motion
$(X_t)_{t\geq 0} = (B_{S_t})_{t\geq 0}$, where $(S_t)_{t\geq 0}$
is an $\alpha/2$-stable process with Laplace exponent $\eta(\lambda) = (2\lambda )^{\alpha / 2}$, i.e. 
$$
 \mathbb{E}\big[e^{-\lambda S_t}\big] = e^{-t ( 2 \lambda)^{\alpha / 2}}, \qquad
\lambda , \ \! t \geq 0, 
$$
see, e.g., Theorem~1.3.23 and pages~55-56 in \cite{applebk2}.
Using the fact that $B_{S_t}/\sqrt{S_t}$ follows the
normal distribution $\mathcal{N}(0,1)$ given $S_t$,
 for $d\geq 1$ and $p \in (0,d)$ we have 
\begin{eqnarray}
\nonumber
\E [|X_t|^{-p}] &=&  \E [|B_{S_t}|^{-p}]
\\
\nonumber
&=&\E\left[ S_t^{-p/2} \E\left[ \frac{S_t^{p/2}}{|B_{S_t}|^p} \ \! \bigg| \ \! S_t\right] \right]
\\
\nonumber
&=& \E\Bigg[ S_t^{-p/2} \int_{\mathbb{S}^{d-1}} \mu_d (d \sigma) \int_0^\infty r^{d-1-p}  \frac{e^{-r^2/2}}{(2\pi)^{d/2} } dr \Bigg]
\\
\nonumber
&=& 2 \frac{ 2^{(d-p-2)/2}}{2^{d/2} \Gamma (d/2)} \Gamma ((d-p)/2)
\E\big[ S_t^{-p/2} \big]
\\
\label{inverse}
&=& \frac{C_{\alpha,d,p}}{t^{p/\alpha}},
\qquad t>0, \quad \alpha \in (1,2),  
\end{eqnarray}
 where $\mu_d$ denotes the surface measure on the
 $d$-dimensional sphere $\mathbb{S}^{d-1}$, 
\begin{equation}
\nonumber 
C_{\alpha,d,p} := 2^{1-p} \frac{\Gamma (p/\alpha) \Gamma ((d-p)/2) }{\alpha \Gamma (p/2)\Gamma (d/2)}, 
\end{equation} 
 and we used the relation
 $\mathbb{E}\big[ S_t^{-p} \big]
 = \alpha^{-1} 2^{1-p} t^{-2p/\alpha} \Gamma (2p/\alpha)  / \Gamma (p)$, $p,t>0$,
see, e.g., Relation~(1.10) in \cite{penent}. 
\subsubsection*{Integration by parts} 
The stochastic representation of
the gradient $\nabla u(x)$ will rely on an 
integration by parts argument. For this,
we will use the weight functions
 $\mathcal{W}_{B(0,R)}(t,x,y)$ and 
 $\mathcal{W}_{\partial B(0,R)} (x,y)$ defined as 
\begin{equation} 
  \label{fjkdl3}
  \mathcal{W}_{B(0,R)}(t,x,y) := \frac{\nabla_x G_R(x,y)}{G_R(x,y)}
  \ \ \mbox{and} \ \  
  \mathcal{W}_{\partial B(0,R)} (x,y) := \frac{\nabla_x P_R(x,y)}{P_R(x,y)},
  \quad x,y\in B(0,R). 
\end{equation}
 
  \begin{lemma}
    \label{2fjlkd}
Let $\alpha \in (1,2)$ and $d\geq 2$. 
\begin{enumerate}[a)]
\item
 Given $\phi$ a bounded measurable function on $\real^d \setminus B(0,R)$, 
 the function 
\begin{equation}
  \label{11} 
  \chi^\phi_1(x) := \mathbb{E}\big[ \phi \big(X_{\tau_R (x),x}\big) \big]
 = \int_{\real^d \setminus B(0,R)} P_R(x,y) \phi (y) dy, 
 \quad x\in \widebar{B}(0,R), 
\end{equation} 
belongs to ${\cal C}^1(B(0,R))\cap {\cal C}^0 ( \widebar{B}(0,R))$,
 with 
   $$
 \nabla \chi^\phi_1(x) 
= \mathbb{E}\left[
   \mathcal{W}_{\partial B(0,R)} \big(x,X_{\tau_R (x),x}\big)
  \phi \big(X_{\tau_R (x),x}\big)\right], \quad x \in B(0,R). 
$$ 
\item 
 Given $h$ a bounded continuous function on $\widebar{B}(0,R)$, 
 the function 
\begin{equation}
\label{11-2} 
\chi^h_2(x) := \mathbb{E}\Bigg[ \int_0^{\tau_R (x)} h(X_{t,x})dt \Bigg]
 = \int_{B(0,R)} G_R(x,y) h(y) dy, 
 \quad x\in \widebar{B}(0,R), 
\end{equation} 
 belongs to ${\cal C}^1(B(0,R))\cap {\cal C}^0 ( \widebar{B}(0,R))$,
 with 
$$
 \nabla \chi^h_2(x) 
  =  \mathbb{E}
   \bigg[
     \int_0^{\tau_R (x)}
     \mathcal{W}_{B(0,R)}(x,X_{t,x}) h(X_{t,x}) dt
     \bigg], \quad x \in B(0,R). 
$$ 
\end{enumerate} 
\end{lemma}
\begin{Proof} 
 $(a)$ Using \eqref{P1} and the
 boundedness of $\phi$ on $\real^d \setminus B(0,R)$, 
 we differentiate \eqref{11}
  under the integral sign, to obtain 
 that $\chi^\phi_1$ is
 in ${\cal C}^1(B(0,R))\cap {\cal C}^0 ( \widebar{B}(0,R))$,
 with 
 $$
 \nabla \chi^\phi_1(x)
= \int_{\real^d \setminus B(0,R)} \nabla_x P_R(x,y) \phi (y) dy
 = \mathbb{E}\left[
  \frac{\nabla_x P_R\big(x,X_{\tau_R (x),x}\big)}{P_R\big(x,X_{\tau_R (x),x}\big)}
  \phi \big(X_{\tau_R (x),x}\big)\right], \ x \in B(0,R). 
$$ 
 $(b)$ Using \eqref{G1}, the condition $d\geq 2$ and the relation
\begin{eqnarray*} 
 \chi^h_2(x) & = & \int_{B(0,R)} G_R(x,y) h(y) dy
 \\
  & = & \int_{B(x,R)} 
 \frac{\kappa_\alpha^d}{|z|^{d-\alpha}}
 \int_0^{r_0(x,z-x)} \frac{t^{\alpha/2-1}}{(1+t)^{d/2}} dt
 \ \! h(z-x) dz,  
 \quad x\in \widebar{B}(0,R), 
\end{eqnarray*} 
 we differentiate \eqref{11}
 under the integral sign 
 and integrate by parts, to obtain 
 $$
 \nabla \chi^h_2(x) = \int_{B(0,R)} \nabla_x G_R(x,y) h(y) dy
 =  \mathbb{E}
   \bigg[
     \int_0^{\tau_R (x)} \frac{\nabla_x G_R(x,X_{t,x})}{G_R(x,X_{t,x})} h(X_{t,x}) dt
     \bigg], 
$$ 
 first for $h$ a ${\cal C}^1$ function
 with compact support in $B(0,R)$,
 then by uniform approximation of 
 $h$ continuous with compact support in $\widebar{B}(0,R)$, 
 and finally by pointwise approximation of $h$ bounded continuous
 on $\widebar{B}(0,R)$, using the bound \eqref{greenbound}.
\end{Proof}
\section{Marked branching process} 
\label{s2}
Let $\rho: \mathbb{R}^+ \rightarrow (0,\infty )$ be a
probability density function on $\real_+$,
and let $(q_{l_0,\ldots ,l_m})_{(l_0,\ldots ,l_m)\in {\cal L}_m}$
be a strictly positive probability mass function on ${\cal L}_m$. 
 We consider 
\begin{itemize}
\item an i.i.d. family $(\tau^{i,j})_{i,j\geq 1}$ of random variables
  with distribution $\rho (t)dt$ on $\real_+$
 and tail distribution function 
 ${\widebar{F}} (t) = \int_t^\infty \rho (ds)ds$, $t\geq 0$, 
\item an i.i.d. family $(I^{i,j})_{i,j\geq 1}$ of discrete
  random variables with distribution 
  $$
  \mathbb{P} ( I^{i,j}=(l_0,\ldots ,l_m) ) = q_{l_0,\ldots ,l_m} >0,
  \qquad (l_0,\ldots , l_m)\in {\cal L}_m,
$$

\item
  an independent family $(X^{(i,j)})_{i,j\geq 1}$
  of symmetric $\alpha$-stable processes. 
\end{itemize}
 In addition, the families of random variables
 $(\tau^{i,j})_{i,j\geq 1}$, $(I^{i,j})_{i,j\geq 1}$ and
 $\big(X^{(i,j)}\big)_{i,j\geq 1}$ are assumed to be mutually independent.

 \medskip

 The probabilistic representation for the solution of \eqref{eq:1}
 uses a branching process started 
 from a particle $x\in B(0,R)$ with label $\widebar{1}=(1)$ and mark
 $i \in \{0, \ldots, m \}$, 
which evolves according to
the process $X_{s,x}^{\widebar{1}} = x + X_s^{(1,1)}$,
$s \in [0,T_{\widebar{1}} ]$, with
$T_{\widebar{1}} = \tau^{1,1} \wedge \tau_R (x) = \min \big(
\tau^{1,1} ,\tau_R (x) \big)$, where in the notation 
\begin{equation}
\nonumber 
  \tau_R (x) := \inf\big\{ t \geq 0,~ x + X_t^{(1,1)} \not\in B(0,R) \big\}, 
\end{equation}
 we omit the reference to the label $(1,1)$. 

 \medskip 

If $\tau^{1,1}<\tau_R (x)$, then the process branches at time $\tau^{1,1}$
into new independent copies of
$(X_t)_{t\geq 0}$, each of them started at
$X_{\tau^{1,1},x}^{\widebar{1}}$, and determined by a random sample
$(l_0,\ldots , l_m)\in {\cal L}_m$ of $I^{1,1}$.
Namely, $|l|:=l_0+\cdots +l_m$ new branches
carrying respectively the marks $i=0,\ldots ,m$
are created with the probability $q_{l_0,\ldots ,l_m}$,
where: 
\begin{enumerate}[a)]
\item the first $l_0$ branches
  carry the mark $0$ and
  are indexed by $(1,1),(1,2),\ldots ,(1,l_0)$,
\item for $i=1,\ldots , m$, the next $l_i$ branches
  carry the mark $i$ and
  are indexed by the labels $(1,l_0+\cdots +l_{i-1}+1),\ldots ,(1,l_0 + \cdots + l_i)$. 
\end{enumerate}
Each new particle then follows independently the above mechanism
in such a way that particles at generation $n\geq 1$ are assigned
a label of the form $\widebar{k} = (1,k_2,\ldots ,k_n) \in \mathbb{N}^n$,
and every branch stops when it leaves the domain $B(0,R)$.

\medskip

 Precisely, the particle with label 
 $\widebar{k} = (1,k_2,\ldots ,k_n) \in \mathbb{N}^n$
 is born at the time $T_{\widebar{k}-}$, where 
 $\widebar{k}- := (1,k_2,\ldots ,k_{n-1})$
 represents the label of its parent,
 and its lifetime $\tau^{n,\pi_n(\widebar{k})}$ is the element of index
 $\pi_n(\widebar{k})$ in the i.i.d. sequence
 $(\tau^{n,j})_{j\geq 1}$, which defines an injection
$$\pi_n:\mathbb{N}^n \to \mathbb{N},
\qquad n\geq 1.
$$ 
 The random evolution of the particle of label $\widebar{k}$
is given by
\begin{equation}
\nonumber 
X_{t,x}^{\widebar{k}} := X^{\widebar{k}-}_{T_{\widebar{k}-},x}+X_{t-T_{\widebar{k}-}}^{n,\pi_n(\widebar{k})},
\qquad t\in [T_{\widebar{k}-},T_{\widebar{k}}],
\end{equation} 
where $T_{\widebar{k}} := T_{\widebar{k}-} + \tau^{n,\pi_n(\widebar{k})} \wedge \tau_R \big( {X^{\widebar{k}-}_{T_{\widebar{k}-},x}}\big)$ and
$$
\tau_R \big( {X^{\widebar{k}-}_{T_{\widebar{k}-},x}}\big):= \inf \big\{ t \geq 0,~ X^{\widebar{k}-}_{T_{\widebar{k}-},x}+X_t^{n,\pi_n(\widebar{k})} \not\in B(0,R) \big\}.
$$
 If $\tau^{n,\pi_n(\widebar{k})} < \tau_R \big( {X^{\widebar{k}-}_{T_{\widebar{k}-},x}}\big)$,
we draw a random sample $(l_0,\ldots ,l_m)$
 of $I_{\widebar{k}} := I^{n,\pi_n(\widebar{k})}$
 with the probability $q_{l_0,\ldots ,l_m}$, 
 and the particle $\widebar{k}$ branches into
 $|I^{n,\pi_n(\widebar{k})}|=l_0+\cdots +l_m$ offsprings,
 indexed by $(1,\ldots ,k_n,j)$, $j=1,\ldots ,|I^{n,\pi_n(\widebar{k})}|$
 and respectively carrying the marks $i=0,\ldots ,m$,
 as in point $(b)$ above.
 Namely, the particles whose index ends with an integer between $1$ and $l_0$
 will carry the mark $0$, and those with index ending with an integer between
 $l_0+\cdots +l_{i-1} +1$ and $l_0+\cdots + l_i$ will carry a mark
 $i\in \{1,\ldots ,m\}$. 
 Finally, the mark of the particle $\widebar{k}$ will be denoted by
 $\theta_{\widebar{k}} \in \{0,\ldots ,m \}$.

 \medskip 

 The set of particles dying inside the ball $B(0,R)$ is denoted by $\mathcal{K}^{\circ}$,
 whereas those dying outside of $B(0,R)$ form a set denoted by $\mathcal{K}^{\partial}$.
 For $n\geq 1$, the set $n$-$th$ generation particles that die
 inside the domain $B(0,R)$ is denoted by $\mathcal{K}_n^\circ$,
 and the set of $n$-th generation particles which die outside of $B(0,R)$
 is denoted by $\mathcal{K}_n^\partial$,
 and we let
 $\mathcal{K}_n = \mathcal{K}_n^\circ \cup \mathcal{K}_n^\partial$. 
\begin{definition}
  We denote by $\mathcal{T}_{x,i}$
  the marked branching process, or random marked tree
  constructed above after starting from the position $x\in \real^d$ and mark 
  $i \in \{0,\ldots ,m\}$ on its first branch. 
\end{definition}
 The tree $\mathcal{T}_{x,0}$
 will be used for the stochastic representation of the solution
 $u(x)$ of the PDE \eqref{eq:1}, while the trees $\mathcal{T}_{x,i}$ 
 will be used for the stochastic representation of  $b_i(x) \cdot \nabla u(x)$,
 $i=1,\ldots , m$. 
 The next table summarizes the notation introduced so far.
  
 \medskip

 \begin{center}
 \begin{tabular}{||l | c||}
 \hline
 Object & Notation \\ [0.5ex]
 \hline\hline
 Initial position & $x$ \\
 \hline
 Tree rooted at $x$ with initial mark
 $\theta_{\widebar{1}} =i$ &  $\mathcal{T}_{x,i}$ \\  \hline
 Particle (or label) of generation $n\geq 1$ & $\widebar{k}=(1,k_2,\ldots ,k_n)$\\
 \hline
 First branching time & $T_{\widebar{1}}$\\
 \hline
 Lifespan of a particle & $T_{\widebar{k}} - T_{\widebar{k}-}$ \\
 \hline
Birth time of the particle $\widebar{k}$ & $T_{\widebar{k}-}$ \\
\hline
Death time of the particle $\widebar{k} \in \mathcal{K}^\circ$& $T_{\widebar{k}} = \T_{\widebar{k}-} + \tau^{n,\pi_n(\widebar{k})}$ \\
\hline
Death time of the particle $\widebar{k} \in \mathcal{K}^\partial$& $T_{\widebar{k}} = T_{\widebar{k}-} + \tau_R \big( {X^{\widebar{k}-}_{T_{\widebar{k}-},x}}\big)$ \\
\hline
Position at birth of the particle $\widebar{k}$ & $X^{\widebar{k}}_{T_{\widebar{k}-},x}$\\
\hline
Position at death of the particle $\widebar{k}$ & $X^{\widebar{k}}_{T_{\widebar{k}},x}$ \\
\hline
Mark of the particle $\widebar{k}$ & $\theta_{\widebar{k}}$\\
\hline
Exit time starting from $x \in B(0,R)$ & $\tau_R (x) := \inf\left\{ t \geq 0,~ x + X_t \not\in B(0,R) \right\}$
\\ 
\hline
\end{tabular}
\end{center}

\medskip
\smallskip

\begin{figure}[H]
\begin{center}
\resizebox{0.55\textwidth}{!}{
\begin{tikzpicture}[scale=0.9,grow=right, sloped][H]
\node[ellipse split,draw,purple,text=black,thick]{time \nodepart{lower} position}
    child {
        node[ellipse split,draw,purple,text=black,thick]{time \nodepart{lower} position}
        child{
        node[ellipse split,draw,thick]{... \nodepart{lower} ...}
        edge from parent
        node[above]{label}
        node[below]{mark}
        }
        child{
        node[ellipse split,draw,thick]{... \nodepart{lower} ...}
        edge from parent
        node[above]{label}
        node[below]{mark}
        } 
        edge from parent
        node[above]{label}
        node[below]{mark}
    }
    child {
        node[ellipse split,draw,blue,text=black,thick]{time \nodepart{lower} position}
        edge from parent
        node[above]{label}
        node[below]{mark}
    };
\end{tikzpicture}
}
\end{center}
\caption{Tree labelling and marking conventions.}
\label{f1} 
\end{figure}
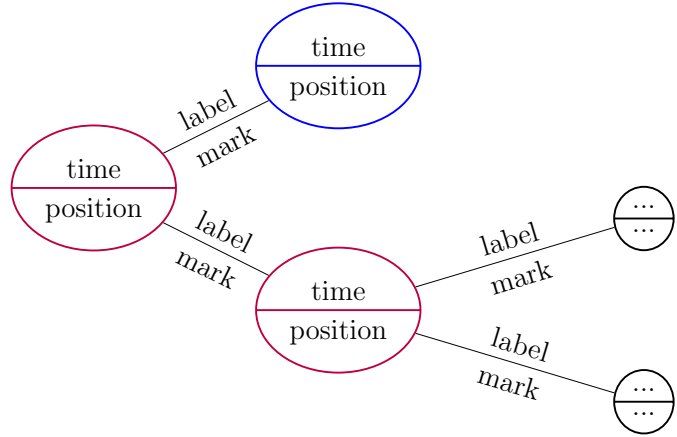

\noindent
Figure~\ref{f1} presents the marking and labeling conventions
used for the graphical representation of random marked trees,
in which different colors represent different ways of branching.  
 
\tikzstyle{level 1}=[level distance=4cm, sibling distance=4cm]
\tikzstyle{level 2}=[level distance=5cm, sibling distance=3cm]

\begin{figure}[H]
\begin{center}
\resizebox{0.85\textwidth}{!}{
\begin{tikzpicture}[scale=0.9,grow=right, sloped]
\node[ellipse split,draw,cyan,thick]{$0$ \nodepart{lower} $x$}
    child {
        node[ellipse split,draw,purple,text=black,thick] {$T_{\widebar{1}}$ \nodepart{lower} $X^{\widebar{1}}_{T_{\widebar{1}},x}$}        
            child {
              node[ellipse split,draw,purple,text=black,thick] {$
                T_{(1,2)}$ \nodepart{lower} $X^{(1,2)}_{T_{(1,2)},x}$} 
                child{
                  node[ellipse split,draw,blue,text=black,thick, right=4cm, below=-1cm]{$
                    T_{(1,2,2)}$ \nodepart{lower} $X^{(1,2,2)}_{T_{(1,2,2)},x}$}
                edge from parent
                node[above]{~~~$(1,2,2)$~~~}
                node[below]{$1$}
                }
                child{
                  node[ellipse split,draw,thick, right=0.3cm]{$
                    T_{(1,2,1)}:=T_{(1,2)}+\tau_R ({X^{(1,2,1)}_{T_{(1,2,1)},x}})$ \nodepart{lower} $X^{(1,2,1)}_{T_{(1,2,1)},x}$}
                edge from parent
                node[above]{$(1,2,1)$}
                node[below]{$0$}
                }
                edge from parent
                node[above] {$(1,2)$}
                node[below]  {$1$}
            }
            child {
              node[ellipse split,draw,blue,text=black,thick,below=-1.6cm] {$
                T_{(1,1)}$ \nodepart{lower} $X^{(1,1)}_{T_{(1,1)},x}$}
                edge from parent
                node[above] {$(1,1)$}
                node[below]  {$0$}
            }
            edge from parent 
            node[above] {$\widebar{1}$}
            node[below]{$0$}
    };
\end{tikzpicture}
}
\end{center}
\caption{Tree labelling and marking conventions.}
\label{f2} 
\end{figure}
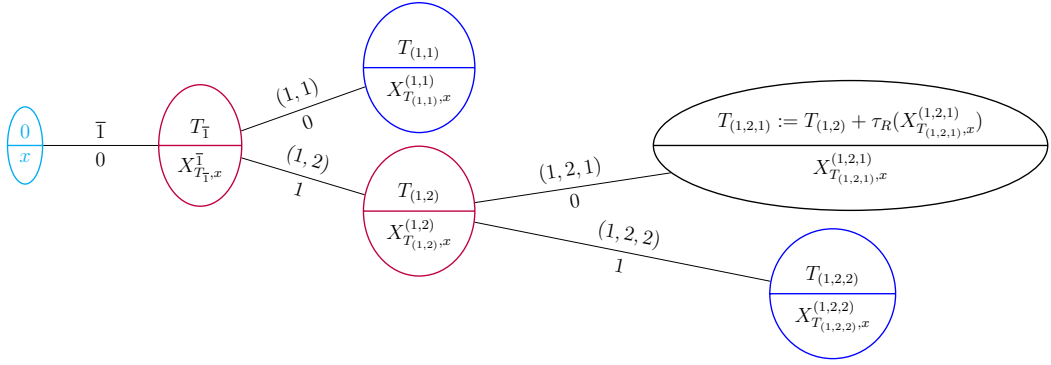

\noindent
 A sample tree for the PDE
$$
\Delta_\alpha  u (t,x) + c_{(0,0)} (x) + c_{(0,1)}(x) u (t,x) \frac{\partial u}{\partial x} (t,x) = 0 $$
in dimension $d=1$ is presented
in Figure~\ref{f2}. Absence of branching is represented in blue,
 branching into two branches,
 one bearing the mark $0$ and the other one bearing the mark $1$, is 
 represented in purple,
 and the black color is used for leaves,
 i.e. for particles that die outside of the domain $B(0,R)$. 

\noindent 
In Figure~\ref{f2} we have 
$\mathcal{K}^{\circ}= \{\widebar{1}, (1,1) , (1,2) ,(1,2,2)\}$
and $\mathcal{K}^{\partial} = \{(1,2,1)\}$.

\section{Probabilistic representation of PDE solutions}
\label{s4} 
\noindent
We consider the weight function $\mathcal{W}(t,x,X)$
defined as 
\begin{equation} 
\label{fjkdl3-2}
\mathcal{W}(t,x,X)
:=
\mathcal{W}_{B(0,R)} (t,x,X_{t,x}) {\bf 1}_{\{X_{t,x} \in B(0,R)\}}
+
\mathcal{W}_{\partial B(0,R)} (x,X_{\tau_R (x),x}) {\bf 1}_{\{X_{t,x} \not\in B(0,R)\}},
\end{equation} 
 $x\in B(0,R)$.
\begin{definition}
  We define the functional $\mathcal{H}_\phi$
  of the random tree $\mathcal{T}_{x,i}$
  with initial mark 
  $\theta_{\widebar{1}} =i \in \{0,\ldots ,m\}$
  as 
\begin{equation}
\label{djsda} 
        { \mathcal{H}_\phi (\mathcal{T}_{x,i}) :=
          \prod_{\widebar{k} \in \mathcal{K}^{\circ}} \frac{c_{I_{\widebar{k}}}\big(X^{\widebar{k}}_{T_{\widebar{k}},x}\big)\mathcal{W}^{\widebar{k}}_{T_{{\widebar{k}}-},x}}{q_{I_{\widebar{k}}}\rho(T_{\widebar{k}} - T_{\widebar{k}-})} \prod_{\widebar{k} \in \mathcal{K}^{\partial}} \frac{   \phi\big(X^{\widebar{k}}_{T_{\widebar{k}},x}\big) \mathcal{W}^{\widebar{k}}_{T_{{\widebar{k}}-},x}}{\widebar{F}(T_{\widebar{k}} - T_{\widebar{k}-})}},
        \quad x\in B(0,R), 
\end{equation} 
 where for
  $\widebar{k} \in \mathcal{K}^{\circ} \cup \mathcal{K}^{\partial}$  
  we let 
\begin{equation}
\label{fjlkds} 
\mathcal{W}^{\widebar{k}}_{T_{{\widebar{k}}-},x}:=
\left\{
\begin{array}{ll} 
  1 & \mbox{if } \theta_{\widebar{k}} = 0,
  \medskip
  \\
  b_{\theta_{\widebar{k}}}
\big(
X^{\widebar{k}}_{T_{{\widebar{k}}-},x}
\big) \cdot \mathcal{W} \big(T_{\widebar{k}}-T_{{\widebar{k}}-},X_{T_{{\widebar{k}}-},x}^{\widebar{k}},X^{\widebar{k}}\big)
& \mbox{if } \theta_{\widebar{k}} = 1,\ldots , m,
\end{array}
\right. 
\end{equation} 
where $\theta_{\widebar{k}} \in \{0,\ldots ,m \}$
denotes the mark of the particle $\widebar{k}$. 
\end{definition} 
\noindent
{\bf Assumption (\hypertarget{B}{$B$})} 
{\em 
 Let $\alpha \in (1,2)$ and $d\geq 2$. 
 We assume that the common probability density function $\rho$
 and tail distribution function ${\widebar{F}}$ 
 of the random times $\tau^{i,j}$'s satisfies the conditions 
\begin{equation} 
\nonumber 
\displaystyle
\sup_{t \in (0,1]} \frac{1}{\rho(t) t^{p/\alpha}} < \infty
\quad
\mbox{and}
\quad 
\E \big[ \big( \widebar{F} ( \tau_R (0) ) \big)^{1-p} \big] < \infty
\end{equation}
for some $p\in (1,d)$.
}

\noindent
 When $\alpha \in (1,2)$ and $R$ is sufficiently small,
 Assumption~(\hyperlink{B}{$B$}) 
 is satisfied by any continuous probability density function
 $\rho (t)$ such that 
 $$ \rho(t) \underset{t\to 0}{\sim} \kappa t^{\delta-1},
 $$ 
 for some $\delta \in (0,1 - p/\alpha ]$ and 
   $\kappa >0$,
   and
   $1 / \widebar{F} ( x) \leq e^{\kappa x}$, $x\geq 0$, 
   for some $\kappa > 0$,
   see, e.g., Lemma~6 in \cite{bogdandirichlet}. 
   This includes for example 
   a gamma distribution with shape parameter
   $\delta \in (0,1 - p/\alpha ]$. 
 The goal of this section is to prove the following result,
 which implies Theorem~\ref{t1.0}. 
\begin{theorem}
\label{t3}
Let $\alpha \in (1,2)$ and $d\geq 2$. 
 Under Assumptions~(\hyperlink{A}{$A$})-(\hyperlink{B}{$B$}), 
 if $R>0$
 and $\max_{l \in {\cal L}_m} \Vert c_l\Vert_\infty $ are 
  sufficiently small, 
  the semilinear elliptic PDE \eqref{eq:1} admits a viscosity solution in
 ${\cal C}^1(B(0,R)) \cap {\cal C}^0 ( \widebar{B}(0,R))$, 
  represented as
  \begin{equation}
    \label{djlkd23}
      u(x) := \mathbb{E} [ \mathcal{H}_\phi (\mathcal{T}_{x,0}) ],
  \qquad x\in B(0,R). 
\end{equation}
\end{theorem}
Before giving the proof of Theorem~\ref{t3} at the end of this section,
we need to state and prove Propositions~\ref{t4.1} and \ref{t4.3}
below. 
 First, in Proposition~\ref{t4.1} we obtain a probabilistic representation for
 the solutions of semilinear elliptic PDEs of the form
 \eqref{eq:1}
 under uniform integrability conditions on
 $(\mathcal{H}_\phi (\mathcal{T}_{x,i}))_{x\in B(0,R)}$,
 $i=0,\ldots , m$.
 Then, in Proposition~\ref{t4.3}
 we show that such conditions are satisfied
 under Assumptions~(\hyperlink{A}{$A$})-(\hyperlink{B}{$B$}). 
\begin{prop} 
\label{t4.1}
Let $\alpha \in (1,2)$ and $d \geq 2$, and 
assume that the family
$(\mathcal{H}(\mathcal{T}_{x,i}))_{x\in B(0,R)}$ is uniformly integrable,
$i=0,\ldots ,m$.
  Then, the function $u(x)$ defined as 
  $$
  u(x) := \mathbb{E} [ \mathcal{H}_\phi (\mathcal{T}_{x,0}) ],
  \quad x\in \widebar{B}(0,R),
  $$
  is a viscosity solution in $C^1(B(0,R)) \cap C^0(\widebar{B}(0,R))$
  of \eqref{eq:1}. In addition, the gradient
  $b_i(x)\cdot \nabla u(x)$ can be represented as the expected value
$$
 b_i(x)\cdot \nabla u(x) = \mathbb{E}\big[ \mathcal{H}_\phi (\mathcal{T}_{x,i})\big], 
 \quad x \in B(0,R), \quad i=1,\ldots ,m. 
 $$ 
\end{prop} 
\begin{Proof}
 Let  
$$
 v_i(x) := \mathbb{E} [ \mathcal{H}_\phi (\mathcal{T}_{x,i}) ], \quad x\in B(0,R),
 \quad i=1,\ldots , m. 
$$
 By considering the first branching at time $T_{\widebar{1}}$
 and letting $\mathcal{T}^{(j)}_{X_{T_{\widebar{1}},x}^{\widebar{1}},i}$,
$j=1+l_0+\cdots +l_{i-1}, \ldots , l_0+\cdots +l_i$,
denote independent tree copies started at
$X_{T_{\widebar{1}},x}^{\widebar{1}}$ with the mark $i\in \{ 0,\ldots , m \}$,
we have 
\begin{align} 
  \nonumber
  & u(x) = 
   \mathbb{E}[\mathcal{H}_\phi (\mathcal{T}_{x,0})]
   \\
  \nonumber
    & =  
  \mathbb{E}\left[
           {\bf 1}_{\{ T_{\widebar{1}} = \tau_R (x) \} }
           \frac{\phi\big(X_{\tau_R (x),x}^{\widebar{1}}\big)}{\widebar{F}( T_{\widebar{1}})}
+ 
    {\bf 1}_{\{ T_{\widebar{1}} < \tau_R (x) \} }
      \sum_{l \in \mathcal{L}_m}
        {\bf 1}_{\{ I_{\widebar{1}} = (l_0,\ldots , l_m) \}}
        \frac{c_{I_{\widebar{1}}}\big(X_{T_{\widebar{1}},x}^{\widebar{1}}\big) }{q_{{I_{\widebar{1}}}} \rho(T_{\widebar{k}})}
    \prod_{i=0}^m
    \prod_{j=1+l_0+\cdots + l_{i-1}}^{l_0+\cdots + l_i} 
    \hskip-0.3cm
    \mathcal{H}_\phi \big(\mathcal{T}^{(j)}_{X_{T_{\widebar{1}},x}^{\widebar{1}},i} \big) 
    \right]
\\
\nonumber 
 & =  
\mathbb{E}\left[
  \phi\big( X_{\tau_R (x),x}^{\widebar{1}} \big) + \int_0^{\tau_R (x)}
   \sum_{l \in \mathcal{L}_m} c_l \big(X_{t,x}^{\widebar{1}}\big) u^{l_0} \big(X_{t,x}^{\widebar{1}}\big) \prod_{i=1}^m v_i^{l_i} \big(X_{t,x}^{\widebar{1}}\big) dt\right]   \\
   \label{fkldsf}
   & = 
   \mathbb{E}\big[
  \phi\big( X_{\tau_R (x),x}^{\widebar{1}} \big)
     \big]
   + \mathbb{E}\bigg[ \int_0^{\tau_R (x)}
    h \big(X_{t,x}^{\widebar{1}}\big) dt\bigg], 
\end{align} 
 where $u(x)$ and the function  
$$
 h(x) := 
 \sum_{l \in \mathcal{L}_m} c_l (x) u^{l_0} (x)
 \prod_{i=1}^m v_i^{l_i} (x), \qquad x\in B(0,R), 
$$
 are bounded continuous on $\widebar{B}(0,R)$ by Lemma~\ref{pl1}. 
 Hence by Lemmas~\ref{2fjlkd} and \ref{pl1} 
 the function $u(x)$ is differentiable in $x\in B(0,R)$,
 with 
\begin{eqnarray*} 
 \nabla u(x) & = &  
 \nabla
 \mathbb{E}\big[
     \phi\big( X_{\tau_R (x),x}^{\widebar{1}} \big)
     \big]
   +
   \nabla
     \mathbb{E} \bigg[ \int_0^{\tau_R (x)}
   h \big(X_{t,x}^{\widebar{1}}\big) dt\bigg]
   \\
   & = &
   \mathbb{E}\big[
  \mathcal{W}_{\partial B(0,R)} \big(x,X_{\tau_R (x),x}^{\widebar{1}}\big)
  \phi \big(X_{\tau_R (x),x}^{\widebar{1}}\big)\big]
   +
   \mathbb{E}
   \bigg[
     \int_0^{\tau_R (x)} \mathcal{W}_{B(0,R)} \big(t,x,X_{t,x}^{\widebar{1}}\big)
     h \big( X_{t,x}^{\widebar{1}} \big)
      dt
      \bigg]
   \\
   & = &
 \mathbb{E}\big[ \mathcal{H}_\phi ( \mathcal{T}_{x,0} ) \mathcal{W}(T_{\widebar{1}},x,X) \big], 
\end{eqnarray*} 
 and by \eqref{fjlkds}-\eqref{djsda} we have 
\begin{align*} 
  b_i(x) \cdot \nabla u(x) 
    & =  \mathbb{E}\big[ \mathcal{H}_\phi ( \mathcal{T}_{x,0} )
     \ b_i(x) \cdot \mathcal{W}(T_{\widebar{1}},x,X) \big]
    \\
    & = 
    \mathbb{E} [ \mathcal{H}_\phi (\mathcal{T}_{x,i}) ]
    \\
    & = 
    v_i(x), \qquad x\in B(0,R),
      \quad i=1,\ldots , m. 
\end{align*} 
 Therefore, using \eqref{jfklds}, Relation~\eqref{fkldsf} rewrites as 
$$
 u(x) = \mathbb{E}\bigg[ \phi \big(X_{\tau_R (x),x}^{\widebar{1}}\big)
   + \int_0^{\tau_R (x)} f\big( X_{t,x}^{\widebar{1}} ,u\big(X_{t,x}^{\widebar{1}}\big),\nabla u\big(X_{t,x}^{\widebar{1}}\big)\big) dt \bigg],
 \quad x\in B(0,R). 
$$
 It then follows from a classical argument
  that $u$ is a viscosity solution of \eqref{eq:1}.
 Indeed, for any $\delta > 0$, by the Markov property we also have
$$
 u(x) = \mathbb{E}\bigg[ u \big(X_{\delta \wedge \tau_R (x),x}^{\widebar{1}}\big)
   + \int_0^{\delta \wedge \tau_R (x)} f\big( X_{t,x}^{\widebar{1}} ,u\big(X_{t,x}^{\widebar{1}}\big),\nabla u\big(X_{t,x}^{\widebar{1}}\big)\big) dt \bigg],
 \quad x\in B(0,R). 
$$
 Next, let $ \xi \in \mathcal{C}^2(B(0,R))$ such that $x$ is a maximum point of $u-\xi$ and $u(x) = \xi(x)$. By the It\^o-Dynkin formula,
 we get 
$$
\mathbb{E}\big[\xi\big(X^{\widebar{1}}_{\delta \wedge \tau_R (x),x}\big)\big]
 = \xi(x)
 + \mathbb{E}\bigg[  \int_0^{\delta \wedge \tau_R (x)} \Delta_\alpha \xi \big(X^{\widebar{1}}_{t,x}\big) dt \bigg]. 
$$
Thus, since $u(x) = \xi(x)$ and $u\leq \xi$, we find 
$$
\mathbb{E}\bigg[ \int_0^{\delta \wedge \tau_R (x)}
  \big( 
  \Delta_\alpha \xi \big(X^{\widebar{1}}_{t,x}\big)
  + f\big(X^{\widebar{1}}_{t,x} , u\big(X^{\widebar{1}}_{t,x}\big),\nabla u\big(X_{t,x}^{\widebar{1}}\big)\big)
  \big) dt \bigg] \geq 0. 
$$
Since $X_{t,x}$ converges in distribution to the constant $x\in \real^d$
as $t$ tends to zero,
it admits an  almost surely convergent subsequence,
hence by continuity and boundedness of $f( \ \! \cdot \ \! ,u( \ \! \cdot  \ \! ))$
together with the mean-value and dominated convergence theorems,
we have 
 $$
 \Delta_\alpha \xi(x)  + f(x,\xi(x),\nabla \xi (x) ) \geq 0, 
$$
 hence $u$ is a viscosity subsolution
 (and similarly a viscosity supersolution)
 of \eqref{eq:1}.
\end{Proof}
 The proof of the next lemma uses the filtration $(\mathcal{F}_n)_{n\geq 1}$
 defined by 
$$
\mathcal{F}_n := \sigma\bigg(T_{\widebar{k}},I_{\widebar{k}},X^{\widebar{k}},\widebar{k} \in \bigcup_{i=1}^n \mathbb{N}^i\bigg), \qquad n \geq 1. 
$$
 Recall that $\mathcal{K}_i^\circ$ (resp. $\mathcal{K}_i^\partial$),
 $i=1,\ldots , n+1$, denotes the set of $i$-th generation particles  
 which die inside (resp. outside) the domain $B(0,R)$,
 and 
 $\mathcal{K}_n = \mathcal{K}_n^\circ \cup \mathcal{K}_n^\partial$. 
\begin{lemma}
\label{ll2}
 Let $p\geq 1$,
 and let $v:B(0,R)\to \real_+$ be a bounded measurable function
 satisfying the inequality 
$$
 v(x) \geq K_1 \E \big[ \big( \widebar{F} ( \tau_R (x) ) \big)^{1-p} \big] 
 + \mathbb{E}\bigg[ 
     \int_0^{\tau_R (x)}
         \big( K_2^p {\bf 1}_{[0,1]}(t) \rho (t) 
         + K_3^p 
         {\bf 1}_{(1,\infty )}(t) 
         \big) 
         \sum_{l = (l_0, \ldots , l_m ) \in \mathcal{L}_m}
         \frac{v^{|l|}(X^{\widebar{1}}_{t,x})}{q_l^{p-1}} dt \bigg],
    $$
  $x\in B(0,R)$, for some
   $K_1, K_2 , K_3 > 0$,
   where $|l| = l_0+\cdots +l_m$. 
 Then, we have 
\begin{equation} 
\label{wh} 
v(x) \geq  \mathbb{E} \left[
  \prod_{\widebar{k} \in \mathcal{K}^{\partial}} \frac{K_1}{\widebar{F}^p(T_{\widebar{k}}-T_{\widebar{k}-})}
  \prod_{\widebar{k} \in \mathcal{K}^{\circ}
    \atop
    T_{\widebar{k}}-T_{\widebar{k}-} \leq 1 }
  \frac{K_2^p}{q^p_{I_{\widebar{k}}} }
    \prod_{\widebar{k} \in \mathcal{K}^{\circ}
    \atop
    T_{\widebar{k}}-T_{\widebar{k}-} > 1 }
  \frac{K_3^p}{q^p_{I_{\widebar{k}}} \rho^p (T_{\widebar{k}}-T_{\widebar{k}-})}
 \right],
\quad
 x\in B(0,R). 
\end{equation} 
\end{lemma} 
\begin{Proof}
 Since $T_{\widebar{1}}$ is independent of 
$\big(X^{\widebar{1}}_{s,x}\big)_{s\geq 0}$ and has the probability density $\rho$, letting
$$
 g (y):= \sum_{l = (l_0, \ldots , l_m ) \in \mathcal{L}_m} \frac{y^{|l|}}{q^{p-1}_l}, 
$$
 we have 
\begin{align}
  \nonumber
 & v(x) \geq 
    \mathbb{E}\bigg[ 
      K_1 \big( \widebar{F} ( \tau_R (x) ) \big)^{1-p} 
       +    \int_0^{\tau_R (x)}
         \big( K_2 {\bf 1}_{[0,1]}(t) \rho (t) 
         + K_3 
         {\bf 1}_{(1,\infty )}(t) 
         \big)^p 
         g \big(v\big(X^{\widebar{1}}_{t,x}\big)\big)
         dt \bigg]
\\
\nonumber
    & = 
    \mathbb{E}\left[
    \mathbb{E}\left[
      K_1 \big( \widebar{F} ( \tau_R (x) ) \big)^{1-p}
      + \int_0^{\tau_R (x)} 
 \big( K_2 {\bf 1}_{[0,1]}(t) \rho (t) 
         +
         K_3 {\bf 1}_{(1,\infty )}(t) 
         \big)^p 
         g \big(v\big(X^{\widebar{1}}_{t,x}\big)\big)
       dt
     \ \! \bigg|
      \ \!
      \big(X^{\widebar{1}}_{s,x}\big)_{s\geq 0}
      \right]\right]
\\
\nonumber
    & = 
  \mathbb{E}\left[
      \mathbb{E}\left[
        \frac{ K_1}{\widebar{F}^p (\tau_R (x))}
                {\bf 1}_{\{ T_{\widebar{1}}= \tau_R (x) \} }
  + 
        \int_0^{\tau_R (x)} 
   \left( K_2^p {\bf 1}_{[0,1]}(t)  
   +
   \frac{K_3^p}{\rho (t) }
    {\bf 1}_{(1,\infty )}(t) 
                  \right) 
              g \big(v\big(X^{\widebar{1}}_{t,x}\big)\big)
          \rho (t) dt
     \ \! \bigg|
      \ \!
      \big(X^{\widebar{1}}_{s,x}\big)_{s\geq 0}
      \right]\right]
 \\
\nonumber
    & =
\mathbb{E}\left[ \frac{ K_1}{\widebar{F}^p(\tau_R (x))}
   {\bf 1}_{\{ T_{\widebar{1}}= \tau_R (x)\} }
  +
   K_2^p g \big(v\big(X^{\widebar{1}}_{T_{\widebar{1}},x}\big)\big)
         {\bf 1}_{\{ T_{\widebar{1}} \leq \min ( 1 , \tau_R (x) ) \} }
           +
  \frac{K_3^p}{\rho (T_{\widebar{1}}) }
   g \big(v\big(X^{\widebar{1}}_{T_{\widebar{1}},x}\big)\big)
   {\bf 1}_{\{ 1 < T_{\widebar{1}} < \tau_R (x) \} } \right]
 \\
\nonumber 
&=   \mathbb{E}\left[  \frac{K_1}{\widebar{F}^p( T_{\widebar{1}})}
  {\bf 1}_{\{ T_{\widebar{1}}= \tau_R (x) \} }
  +
  \frac{1}{q^p_{I_{\widebar{1}}}}
  \left(
  K_2^p 
     {\bf 1}_{ \{ T_{\widebar{1}} \leq \min ( 1 , \tau_R (x) ) \} }
   +  \frac{K_3^p}{\rho (T_{\widebar{1}}) } 
   {\bf 1}_{ \{  1 < T_{\widebar{1}} < \tau_R (x) \} }
   \right)
   v^{|I_{\widebar{1}}|} \big(X^{\widebar{1}}_{T_{\widebar{1}},x}\big)
   \right], 
\end{align} 
 showing that
\begin{equation}
  \label{fkjdslf3} 
v(x) \geq \mathbb{E} \left[
  \prod_{\widebar{k} \in \mathcal{K}_1^{\partial}}\frac{K_1}{\widebar{F}^p(T_{\widebar{k}}-T_{\widebar{k}-})}
  \prod_{\widebar{k} \in \mathcal{K}_1^{\circ}
    \atop
    T_{\widebar{k}}-T_{\widebar{k}-} \leq 1
  }  \frac{K_2^p}{q^p_{I_{\widebar{k}}} }
  \prod_{\widebar{k} \in \mathcal{K}_1^{\circ}
    \atop
      T_{\widebar{k}}-T_{\widebar{k}-} > 1 
  }  \frac{K_3^p }{q^p_{I_{\widebar{k}}} \rho^p (T_{\widebar{k}}-T_{\widebar{k}-})}
  \prod_{\widebar{k} \in \mathcal{K}_2} v\big(X^{\widebar{k}}_{T_{\widebar{k}-},x}\big) \right],
\end{equation} 
 $x\in B(0,R)$.  
 By repeating this argument for the particles in $\widebar{k}\in \mathcal{K}_2$, we find 
\begin{align*} 
 v\big(X^{\widebar{k}}_{T_{\widebar{k}-},x}\big)
& \geq \mathbb{E}\left[ \frac{K_1 }{\widebar{F}^p( T_{\widebar{k}}-T_{\widebar{k}-})}       {\bf 1}_{\{ X^{\widebar{k}}_{T_{\widebar{k}},x} \notin B(0,R)\}}
   \right.
   \\
   & \quad + \frac{1}{q^p_{I_{\widebar{k}}}}
        \left(
   K_2^p {\bf 1}_{ \{
     T_{\widebar{k}}-T_{\widebar{k}-}
     \leq \min ( 1 , \tau_R (x) ) \} }
    + \frac{K_3^p }{\rho^p ( T_{\widebar{k}}-T_{\widebar{k}-}) } 
    {\bf 1}_{ \{ 1 <
      T_{\widebar{k}}-T_{\widebar{k}-}
       < \tau_R (x) \} }
   \right)
   v^{|I_{\widebar{k}}|} \big(X^{\widebar{k}}_{T_{\widebar{k}},x}\big)
   \ \! \bigg| \ \! \mathcal{F}_1 \Bigg]. 
\end{align*} 
 Plugging this expression in \eqref{fkjdslf3} above and using the tower
 property of the conditional expectation, we obtain
 $$
 v(x) \geq \mathbb{E} \left[
   \prod_{\widebar{k} \in \bigcup_{i=1}^2  \mathcal{K}_i^{\partial}}\frac{K_1}{\widebar{F}^p(T_{\widebar{k}}-T_{\widebar{k}-})}
   \prod_{\widebar{k} \in \bigcup_{i=1}^2 \mathcal{K}_i^{\circ}
        \atop
     T_{\widebar{k}}-T_{\widebar{k}-} \leq 1 
   }  \frac{K_2^p}{q^p_{I_{\widebar{k}}} }
   \prod_{\widebar{k} \in \bigcup_{i=1}^2 \mathcal{K}_i^{\circ}
     \atop
     T_{\widebar{k}}-T_{\widebar{k}-} > 1 
   }  \frac{K_3^p}{q^p_{I_{\widebar{k}}} \rho^p (T_{\widebar{k}}-T_{\widebar{k}-})}
   \prod_{\widebar{k} \in \mathcal{K}_4} v \big( X^{\widebar{k}}_{T_{\widebar{k}-},x } \big) \right], 
$$
 and repeating this process inductively leads to
$$
v(x) \geq \mathbb{E} \left[
  \prod_{\widebar{k} \in \bigcup_{i=1}^n  \mathcal{K}_i^{\partial}}\frac{K_1}{\widebar{F}^p(T_{\widebar{k}}-T_{\widebar{k}-})}
  \prod_{\widebar{k} \in \bigcup_{i=1}^n \mathcal{K}_i^{\circ}
        \atop
     T_{\widebar{k}}-T_{\widebar{k}-} \leq 1 
 }  \frac{K_2^p}{q^p_{I_{\widebar{k}}} }
  \prod_{\widebar{k} \in \bigcup_{i=1}^n \mathcal{K}_i^{\circ}
        \atop
     T_{\widebar{k}}-T_{\widebar{k}-} > 1 
 }  \frac{K_3^p}{q^p_{I_{\widebar{k}}} \rho^p (T_{\widebar{k}}-T_{\widebar{k}-})}
  \prod_{\widebar{k} \in \mathcal{K}_{n+1}} v\big(X^{\widebar{k}}_{T_{\widebar{k}-},x}\big) \right], 
$$
$n\geq 1$.
 Using Fatou's lemma as $n$ tends to infinity,
since all particles become eventually extinct with probability one, 
 we obtain \eqref{wh}. 
\end{Proof}
\begin{lemma} 
\label{l4.3}
Let $\alpha \in (1,2)$, $p\in [1,d)$, $d\geq 2$, and 
$$ 
b_{0,\infty} :=  \max_{1 \leq i \leq m} \sup_{x \in B(0,R)} |b_i(x)|, 
\quad
b_{1,\infty} := \max_{1 \leq i \leq m} \sup_{x \in B(0,R)}\frac{|b_i(x)|}{R-|x|}.
$$ 
 Under Assumptions~(\hyperlink{A}{$A$})-(\hyperlink{B}{$B$}),
 we have the bound 
\begin{equation} 
  \label{fjklf3}
 \mathbb{E}\big[ \big|\mathcal{H}_\phi (\mathcal{T}_{x,i})\big|^p \big]
   \leq 
    \mathbb{E} \left[
               \prod_{\widebar{k} \in \mathcal{K}^{\partial}} 
        \frac{K_1}{
          \widebar{F}^p (T_{\widebar{k}}-T_{\widebar{k}-}) }
      \prod_{\widebar{k} \in \mathcal{K}^{\circ}
        \atop
 T_{\widebar{k}}-T_{\widebar{k}-} \leq 1 }
      \frac{ K_4 \max_{l \in {\cal L}_m } \Vert c_l \Vert_\infty^p 
          }{q^p_{I_{\widebar{k}}}}
      \prod_{\widebar{k} \in \mathcal{K}^{\circ}
        \atop   T_{\widebar{k}}-T_{\widebar{k}-} > 1 
      }
      \frac{ K_3 \max_{l \in {\cal L}_m} \Vert c_l \Vert_\infty^p}{
        q^p_{I_{\widebar{k}}} \rho^p (T_{\widebar{k}}-T_{\widebar{k}-})}
        \right], 
        \end{equation}
  $x\in B(0,R)$, $i = 0,\ldots , m$, 
where
\begin{equation}
\label{k1k2} 
K_1 := \Vert \phi\Vert_\infty^p ( 1 + (d+\alpha)^p b^p_{1,\infty} ),
\quad 
 K_3 :=  1 + d^p b^p_{1,\infty} + d^pb^p_{0,\infty} C_{\alpha , d , p}, 
\end{equation} 
 and
 $$ 
K_4 :=
          \sup_{t\in [0, 1]}\frac{1 + d^pb^p_{1,\infty}}{\rho^p (t) } 
     + 
        d^p b^p_{0,\infty} \sup_{t\in [0, 1]}\frac{ C_{\alpha,d,p}}{\rho^p (t) t^{p/\alpha}}.
$$
\end{lemma}
\begin{Proof}
 For $x \in B(0,R)$, let 
\begin{equation} 
\label{product}
 w_i (x) : =
 \mathbb{E}\big[ \big|\mathcal{H}_\phi (\mathcal{T}_{x,i})\big|^p \big]
 = 
 \mathbb{E}_i \left[ \prod_{\widebar{k} \in \mathcal{K}^{\circ}} \frac{\big|c_{I_{\widebar{k}}}\big(X^{\widebar{k}}_{T_{\widebar{k}},x}\big)\big|^p
     \big| \mathcal{W}^{\widebar{k}}_{T_{{\widebar{k}}-},x}\big|^p}{q^p_{I_{\widebar{k}}}\rho^p (T_{\widebar{k}}-T_{\widebar{k}-})} \prod_{\widebar{k} \in \mathcal{K}^{\partial}} \frac{\big| \phi\big(X^{\widebar{k}}_{T_{\widebar{k}},x}\big)\big|^p
     \big|\mathcal{W}^{\widebar{k}}_{T_{{\widebar{k}}-},x} \big|^p}{{\widebar{F}^p}(T_{\widebar{k}}-T_{\widebar{k}-})}\right]
 ,
\end{equation} 
 where $\mathbb{E}_i$ denotes the
 conditional expectation given that the tree 
 $\mathcal{T}_{x,i}$ is started with the mark
 $i\in \{0,\ldots , m\}$. 
 When $\widebar{k} \in \mathcal{K}^\circ$ has mark 
 $\theta_{\widebar{k}}=0$ we have $\mathcal{W}^{\widebar{k}}_{T_{{\widebar{k}}-},x}=1$, 
 whereas when $\theta_{\widebar{k}}\neq 0$, using
 \eqref{greenbound}, \eqref{fjkdl3-2}-\eqref{fjlkds}
 and the Cauchy-Schwarz inequality, we have 
\begin{eqnarray} 
  \nonumber
  \big|\mathcal{W}^{\widebar{k}}_{T_{{\widebar{k}}-},x}\big|
& \leq & \frac{d \big|b_{\theta_k}\big(X^{\widebar{k}}_{T_{\widebar{k}-},x}\big)\big|}{\min
  \big(
  R - \big|
  X^{\widebar{k}}_{T_{\widebar{k}-},x}
  \big| ,
  \big|X^{\widebar{k}}_{T_{\widebar{k}},x}-X^{\widebar{k}}_{T_{\widebar{k}-},x}\big|
  \big)}
\\
\nonumber
  & \leq & d \max \Bigg(
\frac{\big|b_{\theta_k}\big(X^{\widebar{k}}_{T_{\widebar{k}-},x}\big)\big|}{
  R - \big|
  X^{\widebar{k}}_{T_{\widebar{k}-},x}
  \big|
  }
,
\frac{\big|b_{\theta_k}\big(X^{\widebar{k}}_{T_{\widebar{k}-},x}\big)\big|}{\big|X^{\widebar{k}}_{T_{\widebar{k}},x}-X^{\widebar{k}}_{T_{\widebar{k}-},x}\big|}
\Bigg)
\\
\nonumber 
&\leq  & db_{1,\infty} + \frac{db_{0,\infty}}{\big|X^{\widebar{k}}_{T_{\widebar{k}},x}-X^{\widebar{k}}_{T_{\widebar{k}-},x}\big|}. 
\end{eqnarray}
 Similarly, 
 when $\widebar{k} \in \mathcal{K}^\partial $,
 the definition of $\mathcal{W}_{\partial B(0,R)} (x,y)$
 in \eqref{fjkdl3}, together with the 
 bound \eqref{poissonbound} and the Cauchy-Schwarz inequality, imply
 \begin{equation}
   \label{eq3} 
 \big| \mathcal{W}^{\widebar{k}}_{T_{{\widebar{k}}-},x}\big|\leq (d+\alpha) b_{1,\infty}.
\end{equation} 
 Next, by conditional independence
 given $\mathcal{G} := \sigma\big(\tau^{i,j},I^{i,j} \ : \ i,j \geq 1\big)$ 
 of the terms in the product over 
 $\widebar{k} \in \mathcal{K}^{\circ}$ 
 and 
 $\widebar{k} \in \mathcal{K}^{\partial}$, which  
 involve random terms of the form
 $X^{\widebar{k}}_{T_{\widebar{k}},x}-X^{\widebar{k}}_{T_{\widebar{k}-},x}$ 
 given $T_{\widebar{k}}-T_{\widebar{k}-}$, by \eqref{inverse},
 and \eqref{product}-\eqref{eq3}, 
 we have
\begin{align} 
   \nonumber
& w_i(x)
   \\
   \nonumber
   &  \leq  \mathbb{E} \left[ \prod_{\widebar{k} \in \mathcal{K}^{\circ}} \frac{\Vert c_{I_{\widebar{k}}}\Vert_\infty^p}{q^p_{I_{\widebar{k}}} } \E \Bigg[ \frac{2^p}{\rho^p (T_{\widebar{k}}-T_{\widebar{k}-})}
       \Bigg( 1 + d^p b^p_{1,\infty} + \frac{d^p b^p_{0,\infty}}{
         \big|
         X^{\widebar{k}}_{T_{\widebar{k}},x}-X^{\widebar{k}}_{T_{\widebar{k}-},x}
         \big|^p
       } \Bigg) \ \! \bigg| \ \! \mathcal{G}\Bigg] \prod_{\widebar{k} \in \mathcal{K}^{\partial}} \E \Bigg[ \frac{K_1}{{\widebar{F}^p}(T_{\widebar{k}}-T_{\widebar{k}-})} \ \! \bigg| \ \! \mathcal{G}\Bigg] \right]
  \\
 \nonumber
 & =  \mathbb{E} \left[ \prod_{\widebar{k} \in \mathcal{K}^{\circ}}\left(
   \frac{ \Vert c_{I_{\widebar{k}}}\Vert^p_\infty }{q^p_{I_{\widebar{k}}} \rho^p (T_{\widebar{k}}-T_{\widebar{k}-})}
     \Bigg( 1 + d^p b^p_{1,\infty} + 
   \frac{d^p b^p_{0,\infty} C_{\alpha,d,p}}{(T_{\widebar{k}}-T_{\widebar{k}-})^{p/\alpha}} \Bigg) \right)
   \prod_{\widebar{k} \in \mathcal{K}^{\partial}} \frac{K_1}{{\widebar{F}^p}(T_{\widebar{k}}-T_{\widebar{k}-})}  \right]. 
\end{align} 
 Splitting the terms
 in the product over $\widebar{k} \in \mathcal{K}^{\circ}$
between small and
large values of $T_{\widebar{k}}-T_{\widebar{k}-}$, we get 
\begin{align} 
\nonumber
 & w_i (x)
 \\
 & \leq 
    \mathbb{E} \left[ \prod_{\widebar{k} \in \mathcal{K}^{\circ}}
      \frac{ \Vert c_{I_{\widebar{k}}}\Vert^p_\infty }{q^p_{I_{\widebar{k}}}}
    \left(
    \frac{
      K_3
      }{ \rho^p (T_{\widebar{k}}-T_{\widebar{k}-})}
        {\bf 1}_{\{T_{\widebar{k}}-T_{\widebar{k}-} > 1 \}}
        +         K_4 
          {\bf 1}_{\{T_{\widebar{k}}-T_{\widebar{k}-} \leq 1 \}}
        \right)
        \prod_{\widebar{k} \in \mathcal{K}^{\partial}} \frac{K_1}{{\widebar{F}^p}(T_{\widebar{k}}-T_{\widebar{k}-})}  \right]
    , 
        \end{align}
 $x\in B(0,R)$, which yields \eqref{fjklf3}, $i = 0,\ldots , m$.  
\end{Proof}
 Proposition~\ref{t4.3} provides sufficient conditions for 
 the finiteness of the upper bound \eqref{fjklf3}, and for 
 $(\mathcal{H}_\phi (\mathcal{T}_{x,i}))_{x\in B(0,R)}$
 to be bounded in $L^1 (\Omega )$, uniformly in $x\in B(0,R)$,
 $i=0,\ldots ,m$, as required in Proposition~\ref{t4.1}. 
 \begin{prop}
\label{t4.3}
Let $\alpha \in (1,2)$, $p\in [1,d)$, and $d\geq 2$.
Under Assumptions~(\hyperlink{A}{$A$})-(\hyperlink{B}{$B$}), suppose
that the boundary condition $\phi$ is bounded on $\real^d$ 
and that there exists a bounded strictly positive weak solution $v \in H^{\alpha/2}(\real^d) \cap L^\infty(\real^d)$ to the following partial differential inequality: 
\begin{equation}
\label{PDI}
\begin{cases}
\displaystyle
 \Delta_\alpha  v(x) + \widetilde{K}_2 \sum_{l\in \mathcal{L}_m} v^{|l|}(x) \leq 0,
  \quad x\in B(0,R),
  \medskip
  \\
   v(x) \geq \widetilde{K}_1 >0, \quad x\in \real^d \setminus B(0,R),  
\end{cases}
\end{equation}
where $\widetilde{K}_1>K_1
\E \big[ \widebar{F}^{1-p}(\tau_R(0))\big]
$, $\widetilde{K}_2\geq K_3$,
and $K_1,K_3 > 0$ are given by \eqref{k1k2}. 
 Then, for sufficiently small $\max_{l \in {\cal L}_m} \Vert c_l \Vert_\infty $
 we have the bound 
 \begin{equation}
   \label{jfdkldf1} 
   \E [|\mathcal{H}_\phi (\mathcal{T}_{x,i})|^p] \leq v(x) \leq \Vert v \Vert_\infty < \infty, \qquad x\in B(0,R),
   \quad
 i=0,\ldots ,m. 
\end{equation} 
\end{prop}
\begin{Proof}
 We smooth out $v\in H^{\alpha/2}(\real^d)$ 
 as 
\begin{equation}
\nonumber 
  v_\varepsilon(x) :=
        \frac{1}{\varepsilon}
        \int_{-\infty}^\infty 
        \psi \left( \frac{x-y}{\varepsilon} \right)v(y)dy, \quad x \in \real,
        \quad \varepsilon >0, 
        \end{equation}
where $\psi : \real \to \real_+$ is a mollifier 
such that $\int_{-\infty}^\infty \psi (y) dy =1$.
By \eqref{PDI} and Jensen's inequality, we have 
\begin{eqnarray*}
  \lefteqn{
    \Delta_\alpha  v_\varepsilon (x)
  + \widetilde{K}_2 \sum_{l\in \mathcal{L}_m} v_\varepsilon^{|l|}(x) 
  }
  \\
   & = & 
        \frac{1}{\varepsilon}
        \int_{-\infty}^\infty
        \Delta_\alpha  \psi \left( \frac{x-y}{\varepsilon} \right)v(y)dy
        + \widetilde{K}_2 \sum_{l\in \mathcal{L}_m}
        \left(
                \frac{1}{\varepsilon}
                \int_{-\infty}^\infty
                \psi \left( \frac{x-y}{\varepsilon} \right)v(y)dy
        \right)^{|l|}
  \\
   & \leq & 
        \frac{1}{\varepsilon}
        \int_{-\infty}^\infty 
      \psi \left( \frac{x-y}{\varepsilon} \right)\Delta_\alpha v(y)dy
        + \widetilde{K}_2 \sum_{l\in \mathcal{L}_m}
                \frac{1}{\varepsilon}
        \int_{-\infty}^\infty 
        \psi \left( \frac{x-y}{\varepsilon} \right)v^{|l|}(y)dy
  \\
  & \leq & 0,
  \quad x\in B(0,R). 
\end{eqnarray*} 
Applying the It\^o-Dynkin formula to $v_\varepsilon (X_{s,x})$ with
$ v_\varepsilon \in H^{\alpha}(\real^d)$, by \eqref{PDI} we have 
\begin{eqnarray*} 
 v_\varepsilon(x)
 & = & \mathbb{E}\left[  v_\varepsilon\big( X_{\tau_R (x)}^x \big) -
    \int_0^{\tau_R (x)}
    \Delta_\alpha  v_\varepsilon(X_{t,x}) dt \right]
  \\
  & \geq & \mathbb{E}\left[ \widetilde{K}_1 + \int_0^{\tau_R (x)} \widetilde{K}_2 \sum_{l\in \mathcal{L}_m}
    \frac{v_\varepsilon^{|l|}(X_{t,x})}{q_l^{p-1}} dt \right],
  \quad x\in B(0,R). 
\end{eqnarray*} 
Thus, passing to the limit as $\varepsilon$ tends to zero, 
by dominated convergence and the facts that
 $\E [ \tau_R (x) ] < \infty$ and 
 $v(x)$ is upper and lower bounded in
$(0,\infty)$, for some sufficiently small $K_2>0$ we have   
\begin{align*} 
  v(x) & \geq \widetilde{K}_1 + \mathbb{E}\bigg[ \int_0^{\tau_R (x)} \widetilde{K}_2 \sum_{l\in \mathcal{L}_m} \frac{v^{|l|}(X_{t,x})}{q_l^{p-1}} dt \bigg]
  \\
  & \geq K_1 + \mathbb{E}\bigg[ \int_0^{\tau_R (x)} 
             \big( K_2^p {\bf 1}_{[0,1]}(t) \rho (t)
         + \widetilde{K}_2^p  
         {\bf 1}_{(1,\infty )}(t) 
         \big) 
         \sum_{l\in \mathcal{L}_m} \frac{v^{|l|}(X_{t,x})}{q_l^{p-1}} dt \bigg],
  \quad x\in B(0,R). 
\end{align*} 
 Hence by Lemmas~\ref{ll2} and \ref{l4.3} we have, 
 provided that $\max_{l \in {\cal L}_m} \Vert c_l \Vert_\infty$
 is sufficiently small,
\begin{align*} 
v(x) & \geq  \mathbb{E} \left[
  \prod_{\widebar{k} \in \mathcal{K}^{\partial}} \frac{K_1}{\widebar{F}^p (T_{\widebar{k}}-T_{\widebar{k}-})}
  \prod_{\widebar{k} \in \mathcal{K}^{\circ}
    \atop
    T_{\widebar{k}}-T_{\widebar{k}-} \leq 1 }
  \frac{K_2^p 
  }{q^p_{I_{\widebar{k}}} }
  \prod_{\widebar{k} \in \mathcal{K}^{\circ}
    \atop
    T_{\widebar{k}}-T_{\widebar{k}-} > 1 }
  \frac{\widetilde{K}_2 }{q^p_{I_{\widebar{k}}} \rho^p (T_{\widebar{k}}-T_{\widebar{k}-})}
  \right]
\\
& \geq  \mathbb{E} \left[
  \prod_{\widebar{k} \in \mathcal{K}^{\partial}} \frac{K_1}{\widebar{F}^p(T_{\widebar{k}}-T_{\widebar{k}-})}
  \prod_{\widebar{k} \in \mathcal{K}^{\circ}
    \atop
    T_{\widebar{k}}-T_{\widebar{k}-} \leq 1 }
  \frac{K_4 \max_{l \in {\cal L}_m} \Vert c_l \Vert_\infty^p 
}{q^p_{I_{\widebar{k}}}
  }
  \prod_{\widebar{k} \in \mathcal{K}^{\circ}
    \atop
    T_{\widebar{k}}-T_{\widebar{k}-} > 1 }
  \frac{K_3 \max_{l \in {\cal L}_m} \Vert c_l \Vert_\infty^p 
  }{q^p_{I_{\widebar{k}}} \rho^p (T_{\widebar{k}}-T_{\widebar{k}-})}
   \right]
\\
 & \geq  
\mathbb{E}\big[ \big|\mathcal{H}_\phi (\mathcal{T}_{x,i})\big|^p \big],
\end{align*} 
$x\in B(0,R)$, $i = 0,\ldots , m$,
which yields \eqref{jfdkldf1} from \eqref{product}. 
\end{Proof}

\begin{Proofy}~\ref{t3}.
  By Theorem~1.2 in \cite{penent2}, 
  the partial differential inequality \eqref{PDI} admits a
  positive (continuous) viscosity solution $v(x)$ on $\real^d$
  when $R>0$ is sufficiently small. 
  In addition, by Proposition~3.5 in \cite{penent2}, 
  $v \in H^{\alpha/2}(\real^d) \cap L^\infty(\real^d)$ and is a  
  weak solution of \eqref{PDI}.
  We conclude the proof by the application
  of Propositions~\ref{t4.1} and \ref{t4.3}. 
\end{Proofy}
\noindent 
Lemma~\ref{pl1}
 extends Lemma~3.3 in
 \cite{penent2} from $i=0$ to $i\in \{1,\ldots , m\}$.
\begin{lemma}
\label{pl1}
Let $i\in \{0,\ldots ,m\}$,
and assume  that   
 $(\mathcal{H}(\mathcal{T}_{x,i}))_{x\in B(0,R)}$ is uniformly integrable. 
 Then, the function 
 $v_i(x):= \mathbb{E} [ \mathcal{H}_\phi (\mathcal{T}_{x,i}) ]$
 is continuous in $x\in \widebar{B}(0,R)$.
 \end{lemma} 
\begin{Proofz} {\em (given for completeness)}. 
  Let $x\in \widebar{B}(0,R)$. By Lemma~3.2 therein, 
 for any sequence $(x_n)_{n\in \inte}$ in $B(0,R)$
 converging fast enough to $x \in \widebar{B}(0,R)$ we have
 \begin{eqnarray}
\label{eq:tau}
\mathbb{P}\Big(\lim_{n\rightarrow \infty} \tau_R (x_n) = \tau_R (x) \Big) = 1, 
\end{eqnarray}
 and letting $\tau_{{\widebar{k}},x} := \tau_R \big( {X^{\widebar{k}-}_{T_{\widebar{k}-},x}}\big)$,
 $\widebar{k} \in \mathcal{K}$, the event  
$$
 A_{\widebar{k}} := 
 \left\{ \lim_{n \rightarrow \infty} \tau_{{\widebar{k}},x_n}  = \tau_{{\widebar{k}},x} \right\} \bigcap \left\{  \lim_{n \rightarrow \infty}X_{\cdot, x_n}^{\widebar{k}}= X_{\cdot, x}^{\widebar{k}} \right\}, 
$$
 has probability one. 
 Again, by Lemma~3.2-$a)$ in {\em ibid}, 
 for some $n_0 (\omega )$ large enough we have 
$$
X_{\tau_{{\widebar{k}},x_n} }^{\widebar{k}}
=
X_{\tau_{{\widebar{k}},x} }^{\widebar{k}}
+x_n-x,  
$$
 and $\tau_{{\widebar{k}},x_n}=\tau_{{\widebar{k}},x}$, $n \geq n_0 (\omega )$. 
 Therefore, using the continuity of the functions
 $\phi$ and $c_l,l\in \mathcal{L}$, we have 
$$
 \lim_{n \rightarrow \infty}
 \phi\big(X_{\tau_{{\widebar{k}},x_n} }^{\widebar{k}}\big)
 \mathcal{W}^{\widebar{k}}_{T_{{\widebar{k}}-},x_n}
 \mathbbm{1}_{ \{ T_{\widebar{k}} = \tau_{{\widebar{k}},x_n} \} }
 =
 \phi\big(X_{\tau_{{\widebar{k}},x}}^{\widebar{k}}\big)
 \mathcal{W}^{\widebar{k}}_{T_{{\widebar{k}}-},x}
 \mathbbm{1}_{\{ T_{\widebar{k}} = \tau_{{\widebar{k}},x}\}}, ~~ \mathbb{P}- \text{a.s.}
$$
 and
 $$
 \lim_{n \rightarrow \infty}
 c_{I_{\widebar{k}}}\big(X_{T_{\widebar{k}},x_n}^{\widebar{k}}\big)
 \mathcal{W}^{\widebar{k}}_{T_{{\widebar{k}}-},x_n}
 \mathbbm{1}_{ \{ T_{\widebar{k}} < \tau_{{\widebar{k}},x_n} \} } =  \frac{c_{I_{\widebar{k}}}\big(X_{T_{\widebar{k}},x}^{\widebar{k}}\big)}{q_{I_{\widebar{k}}}}
 \mathcal{W}^{\widebar{k}}_{T_{{\widebar{k}}-},x}
 \mathbbm{1}_{\{ T_{\widebar{k}} < \tau_{{\widebar{k}},x} \}}, ~~ \mathbb{P}- \text{a.s.}. 
$$
 Hence by \eqref{djsda}, 
 on the event $A := \bigcap_{\widebar{k} \in \mathcal{K} } A_{\widebar{k}}$ 
 of probability one, we have 
$$
 \lim_{n \rightarrow \infty}
 \mathcal{H}_\phi (\mathcal{T}_{x_n,i} (\omega))
 =
 \mathcal{H}_\phi (\mathcal{T}_{x,i} (\omega))
. 
$$
Therefore,
 for any sequence $(x_n)_{n\geq 1}$ converging to
  $x\in \widebar{B}(0,R)$ fast enough, we have
$$
 \mathbb{P}\big(\lim_{n \rightarrow \infty}
 \mathcal{H}_\phi (\mathcal{T}_{x_n,i} )
 =
 \mathcal{H}_\phi (\mathcal{T}_{x,i} (\omega) \big) = 1,
$$ 
 which yields $\lim_{n\to \infty} v_i (x_n) =  v_i (x)$ by uniform integrability of
 $(\mathcal{H}_\phi (\mathcal{T}_{x,i} (\omega)))_{x\in B(0,R)}$. 
\end{Proofz}

\tikzstyle{level 1}=[level distance=4cm, sibling distance=3cm]
\tikzstyle{level 2}=[level distance=5cm, sibling distance=5cm]

\section{Numerical examples} 
\label{s5}
In this section,
we consider numerical applications of the probabilistic representation
\eqref{djlkd23}.
The paths of the $\alpha$-stable process $X_t = B_{S_t}$ are
simulated by time discretization,
 by generating independent random samples of
 Brownian motion and of the $\alpha/2$-stable process
 $(S_t)_{t\in \real_+}$ using the identity in distribution  
$$
S_t \simeq 2t^{2/\alpha }\frac{\sin( \alpha \big(U+\pi / 2) /2 )}{\cos^{2/\alpha } (U)} \left(
\frac{\cos (U- \alpha (U+\pi / 2) /2 )}{E}\right)^{-1+2/\alpha }
$$
based on the Chambers-Mallows-Stuck (CMS) method, 
where $U$ is uniform on $(-\pi / 2,\pi /2)$,
and $E$ is exponential with unit parameter, 
see Relation~(3.2) in \cite{weron1996chambers}. 
In order to keep computation times to a reasonable level,
 the probability density $\rho (t)$
of $\tau^{i,j}$, $i,j\geq 1$,
is taken to be gamma with shape parameters ranging from $1.5$ to $1.7$.
 The C codes used to plot 
 Figures~\ref{fig3}-\ref{fig4} 
 are available at 
 {\url{https://github.com/nprivaul/fractional_elliptic}.}

\medskip

 Given $k\geq 0$, we consider the function
$$
\Phi_{k,\alpha} (x) := (1-|x|^2)^{k+\alpha / 2}_+,
\qquad x\in \real^d,
$$
 which is Lipschitz if $k>1-\alpha/2$, and solves the Poisson problem
 $\Delta_\alpha \Phi_{k,\alpha} = -\Psi_{k,\alpha}$
 on $\real^d$, with
 \begin{align}
\nonumber
&  \Psi_{k,\alpha} (x)
\\
\nonumber 
& := \left\{
\begin{array}{ll}
  \displaystyle
 \frac{\Gamma ( ( d+\alpha ) / 2)
   \Gamma (k+1+ \alpha / 2 )}{
   2^{-\alpha}
   \Gamma (k+1)\Gamma ( d / 2 )}
~{}_2F_1\left(
  \frac{d+\alpha}{2},-k;\frac{d}{2};|x|^2
  \right), ~~ |x|\leq 1
  \medskip
  \\
\displaystyle
 \frac{2^{\alpha}\Gamma ( ( d+\alpha ) /2 )
  \Gamma (k+1+ \alpha / 2 )}{\Gamma (k+1+ ( d+\alpha ) / 2 )
   \Gamma (- \alpha / 2 )
  |x|^{d+\alpha}
 }
  {}_2F_1\left(
  \frac{d+\alpha}{2},\frac{2+\alpha}{2};k+1+\frac{d+\alpha}{2};
  \frac{1}{|x|^2}
\right), ~~| x|>1
\end{array}
\right.
\end{align}
 $x\in \real^d$, where ${}_2F_1 ( a,b;c;y)$ is
 Gauss's hypergeometric function,
 see (5.2) in \cite{getoor},
Lemma~4.1 in \cite{biler2015nonlocal},
and Relation~(36) in \cite{oberman}.
\subsubsection*{Linear gradient term} 
\begin{figure}[H]
\resizebox{0.9\textwidth}{!}{
\begin{tikzpicture}[scale=0.9,grow=right, sloped]
\node[ellipse split,draw,cyan,text=black,thick]{$0$ \nodepart{lower} $x$}
    child {
        node[ellipse split,draw,purple,text=black,thick] {$T_{\widebar{1}}$ \nodepart{lower} $X^{\widebar{1}}_{T_{\widebar{1}},x}$}   
            child {
                node[ellipse split,draw,purple,text=black,thick, right=0.cm] {$T_{(1,1)}$ \nodepart{lower} $X^{(1,1)}_{T_{(1,1)},x}$}
                child{
                    node[ellipse split,draw,black, right=0.cm] {$T_{(1,1,1)} : = T_{(1,1)}+\tau_R (X^{(1,1,1)}_{T_{(1,1,1)},x})$ \nodepart{lower} $X^{(1,1,1)}_{T_{(1,1,1)},x}$}
                    edge from parent
                    node[above]{~$(1,1,1)$~}
                    node[below]{1}
                }
                edge from parent
                node[above] {$(1,1)$}
                    node[below]{1}
            }
            edge from parent 
            node[above] {$\widebar{1}$}
                    node[below]{0}
    };
\end{tikzpicture}
} 

\resizebox{0.506\textwidth}{!}{
\begin{tikzpicture}[scale=0.9,grow=right, sloped]
\node[ellipse split,draw,cyan,text=black,thick]{$0$ \nodepart{lower} $x$}
    child {
        node[ellipse split,draw,purple,text=black,thick] {$T_{\widebar{1}}$ \nodepart{lower} $X^{\widebar{1}}_{T_{\widebar{1}},x}$}   
            child {
                node[ellipse split,draw,blue,text=black,thick, right=0.cm] {$T_{(1,1)}$ \nodepart{lower} $X^{(1,1)}_{T_{(1,1)},x}$}
                edge from parent
                node[above] {$(1,1)$}
                node[below]{1}
            }
            edge from parent 
            node[above] {$\widebar{1}$}
            node[below]{0}
    };
\end{tikzpicture}
} 
\caption{Random tree samples for the PDE \eqref{eq:lg}.} 
\label{f1-2} 
\end{figure}
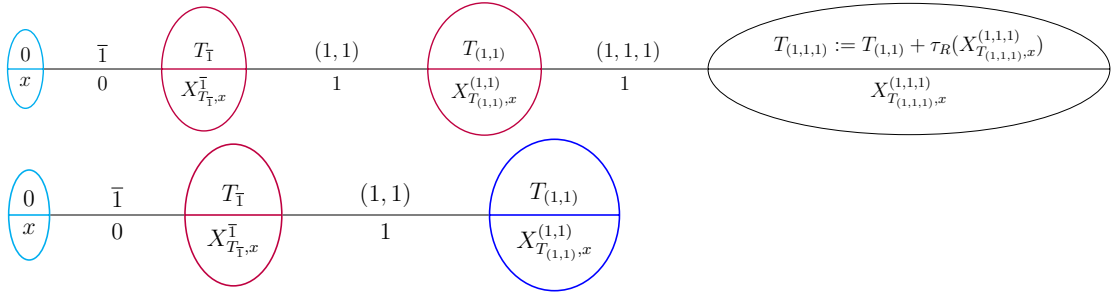

 We take $R=1$, $m=1$, ${\cal L}_1 = \{ (0,0) , (0,1) \}$, and 
 $$
 c_{(0,0)}(x) := \Psi_{k,\alpha}(x) + (2k+\alpha)|x|^2(1-|x|^2)^{k+\alpha/2},
 \quad
 c_{(0,1)}(x) := 1, 
 \quad
 b_1(x) := (1-|x|^2) x, 
$$
 and consider the PDE
\begin{equation}
  \label{eq:lg}
 \Delta_\alpha u (x) + \Psi_{k,\alpha}(x) + (2k+\alpha)|x|^2(1-|x|^2)^{k+\alpha/2} + (1-|x|^2) x \cdot \nabla u(x)  = 0, 
\end{equation}
 $x \in B(0,1)$, with $u(x) =0$ for $x \in \real^d \setminus B(0,1)$,
 and explicit solution 
\begin{equation}
\nonumber 
   u(x) = \Phi_{k,\alpha} (x) 
   = (1-| x|^2)^{k+\alpha/2}_+, \qquad
x \in \real^d. 
\end{equation} 
\noindent  
The random tree associated to \eqref{eq:lg} starts at the point $x\in B(0,1)$, 
and branches into \textcolor{blue}{0 branch} or \textcolor{purple}{1 branch} 
as in the random tree samples of Figure~\ref{f1-2}. 

\medskip

\vskip-0.2cm

\begin{figure}[H]
\centering
\hskip-0.4cm
\begin{subfigure}{.5\textwidth}
\centering
\includegraphics[width=\textwidth]{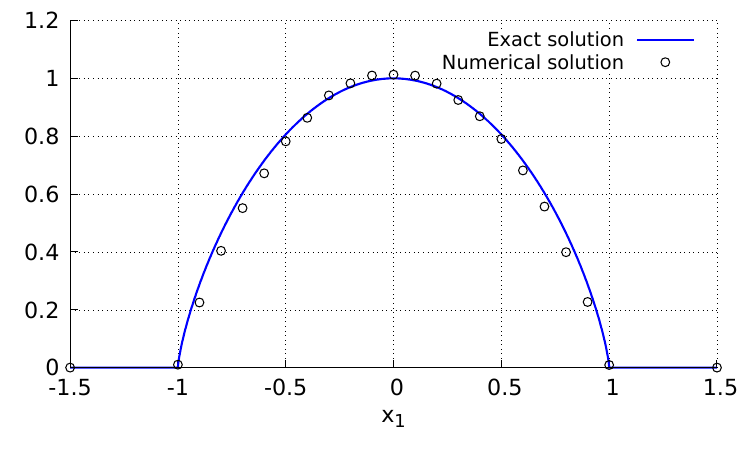}
\vskip-0.3cm
\caption{Numerical solution of \eqref{eq:lg} with $k=0$.}
\label{3a}
\end{subfigure}
\begin{subfigure}{.50\textwidth}
\centering
\includegraphics[width=\textwidth]{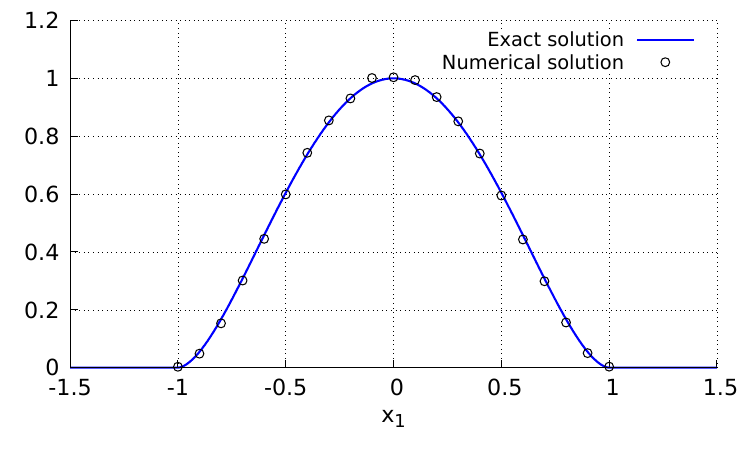}
\vskip-0.3cm
\caption{Numerical solution of \eqref{eq:lg} with $k=1$.}
\label{3b}
\end{subfigure}
\caption{Numerical solution of \eqref{eq:lg} in dimension $d=10$ with $\alpha =1.75$.}
\label{fig3}
\end{figure}

\noindent
 The simulations of Figures~\ref{fig3}-\subref{3a}) and \ref{fig3}-\subref{3b})  respectively use $10^7$ and $2\times 10^7$ Monte Carlo samples. 

\medskip

\subsubsection*{Nonlinear gradient term} 
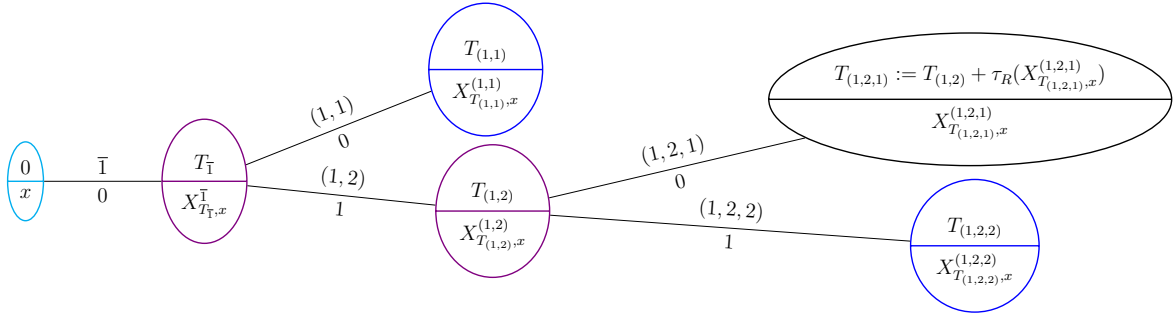
\begin{figure}[H]
\begin{center}
\resizebox{0.95\textwidth}{!}{
\begin{tikzpicture}[scale=0.9,grow=right, sloped]
\node[ellipse split,draw,cyan,text=black,thick]{$0$ \nodepart{lower} $x$}
    child {
        node[ellipse split,draw,violet,text=black,thick] {$T_{\widebar{1}}$ \nodepart{lower} $X^{\widebar{1}}_{T_{\widebar{1}},x}$}   
            child {
                node[ellipse split,draw,violet,text=black,thick,right=1.3cm,below=-3cm] {$T_{(1,2)}$ \nodepart{lower} $X^{(1,2)}_{T_{(1,2)},x}$} 
                child{
                node[ellipse split,draw,blue,text=black,thick,right=5.2cm,below=-2.9cm]{$T_{(1,2,2)}$ \nodepart{lower} $X^{(1,2,2)}_{T_{(1,2,2)},x}$}
                edge from parent
                node[above]{$(1,2,2)$}
                node[below]{1}
                }
                child{
                node[ellipse split,draw,thick, right=1.04cm]{$T_{(1,2,1)}:=T_{(1,2)}+ \tau_R (X^{(1,2,1)}_{T_{(1,2,1)},x})$ \nodepart{lower} $X^{(1,2,1)}_{T_{(1,2,1)},x}$}
                edge from parent
                node[above]{$(1,2,1)$}
                node[below]{0}
                }
                edge from parent
                node[above] {$(1,2)$}
                node[below]{1}
            }
            child {
                node[ellipse split,draw,blue,text=black,thick, right=0.cm] {$T_{(1,1)}$ \nodepart{lower} $X^{(1,1)}_{T_{(1,1)},x}$}
                edge from parent
                node[above] {$(1,1)$}
                node[below]{0}
            }
            edge from parent 
            node[above] {$\widebar{1}$}
                node[below]{0}
    };
\end{tikzpicture}
} 
\end{center}
\caption{Random tree sample for the PDE \eqref{eq:nlg}.} 
\label{f1-3} 
\end{figure}

 In this example we take ${\cal L}_1 = \{ (0,0), (0,2) \}$,  
 $$
 c_{(0,0)}(x) := \Psi_{k,\alpha}(x) + (2k+\alpha)^2|x|^4(1-|x|^2)^{2k+\alpha}, 
 \quad
 c_{(0,2)}(x) := - 1, 
 \quad
 b_1(x) := (1-|x|^2) x, 
$$
 and consider the PDE with nonlinear gradient term
\begin{equation}
\label{eq:nlg}
 \Delta_\alpha u (x) + \Psi_{k,\alpha}(x) + (2k+\alpha)^2|x|^4(1-|x|^2)^{2k+\alpha} - ((1-|x|^2)  x\cdot \nabla u(x) )^2 = 0,
\end{equation}
$x \in B(0,1)$,
with $u(x) =0$ for $x \in \real^d \setminus B(0,R)$,
and explicit solution
$$
u(x) = \Phi_{k,\alpha} (x) 
= (1-| x|^2)^{k+\alpha/2}_+,
\qquad x \in \real^d.
$$ 
The random tree associated to \eqref{eq:nlg} starts at a point $x\in B(0,1)$
and branches into \textcolor{blue}{0 branch} 
 or \textcolor{violet}{2 branches} 
 as in the random tree sample of Figure~\ref{f1-3}. 

\medskip

\vskip-0.3cm

\begin{figure}[H]
\centering
\hskip-0.4cm
\begin{subfigure}{.5\textwidth}
\centering
\includegraphics[width=\textwidth]{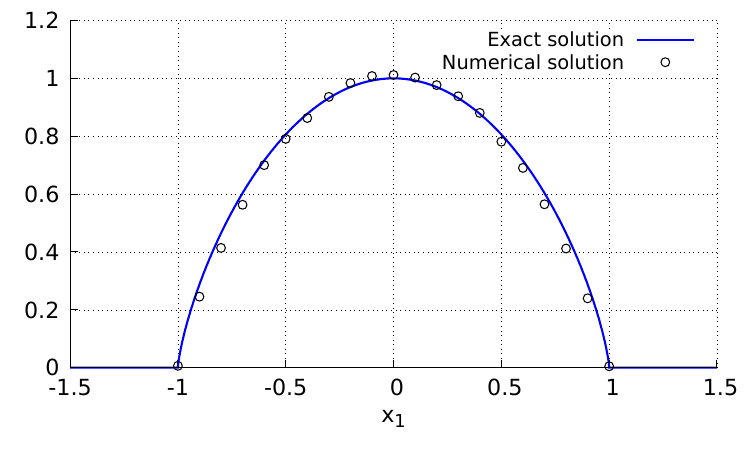}
\vskip-0.3cm
\caption{Numerical solution of \eqref{eq:nlg} with $k=0$.}
\label{4a}
\end{subfigure}
\begin{subfigure}{.50\textwidth}
\centering
\includegraphics[width=\textwidth]{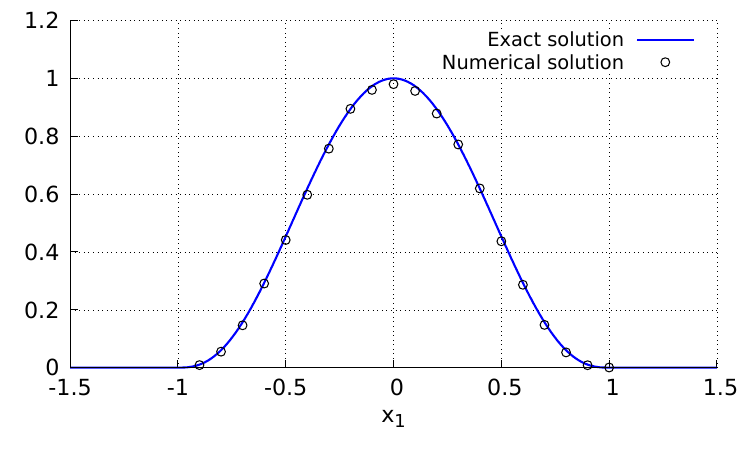}
\vskip-0.3cm
\caption{Numerical solution of \eqref{eq:nlg} with $k=2$.}
\label{4b}
\end{subfigure}
\caption{Numerical solution of \eqref{eq:nlg} in dimension $d=10$ with $\alpha =1.75$.}
\label{fig4}
\end{figure}

\noindent
 The simulations of Figure~\ref{fig4} use five million Monte Carlo samples. 

\footnotesize

\def\cprime{$'$} \def\polhk#1{\setbox0=\hbox{#1}{\ooalign{\hidewidth
  \lower1.5ex\hbox{`}\hidewidth\crcr\unhbox0}}}
  \def\polhk#1{\setbox0=\hbox{#1}{\ooalign{\hidewidth
  \lower1.5ex\hbox{`}\hidewidth\crcr\unhbox0}}} \def\cprime{$'$}

\end{document}